\documentclass[12pt]{article}
\usepackage{fullpage}

\usepackage{amsmath,amssymb}
\usepackage{graphicx}
\usepackage{color}
\usepackage{cite}
\usepackage{curves}
\usepackage{subfigure}

\newcommand{\ve}[1]{\text{\boldmath $#1$}}
\newcommand{\sd}{\mathrm{d}}
\newcommand{\ii}{\mathrm{i}}
\newcommand{\pd}[2]{\frac{\partial #1}{\partial #2}}
\newcommand{\pdd}[2]{\frac{\partial^2 #1}{\partial #2^2}}

\newcommand{\od}[2]{\frac{\sd #1}{\sd #2}}
\newcommand{\ep}{\epsilon}

\DeclareMathOperator{\sn}{sn}
\DeclareMathOperator{\dn}{dn}
\DeclareMathOperator{\K}{K}

\DeclareMathOperator{\im}{Im}
\DeclareMathOperator{\Cc}{Cc}
\DeclareMathOperator{\sech}{sech}

\newcommand{\f}[2]{\frac{#1}{#2}}

\newcommand{\phit}{\tilde{\phi}}

\newcommand{\xt}{\tilde{x}}
\newcommand{\yt}{\tilde{y}}
\newcommand{\nt}{\tilde{n}}
\newcommand{\zt}{\tilde{z}}
\newcommand{\ft}{\tilde{f}}
\newcommand{\zetat}{\tilde{\zeta}}
\newcommand{\xit}{\tilde{\xi}}
\newcommand{\etat}{\tilde{\eta}}

\newcommand{\nuh}{{\hat{\nu}}}
\newcommand{\Yh}{{\hat{Y}}}
\newcommand{\Jh}{{\hat{J}}}
\newcommand{\jod}{{j_{\text{1D}}}}

\title{Two-dimensional modelling of electron flow through a poorly conducting layer}
\date{}
\author{J. P. Black\thanks{black@maths.ox.ac.uk}
\quad\quad
C. J. W. Breward
\quad\quad
P. D. Howell\\
Mathematical Institute, University of Oxford,\\
Andrew Wiles Building, Oxford OX2 6GG, UK}%\\

\begin{document}
\maketitle
%\slugger{siap}{xxxx}{xx}{x}{x--x}%slugger should be set to mms, siap, sicomp, sicon, sidma, sima, simax, sinum, siopt, sisc, or sirev

\begin{abstract}
Motivated by contact resistance on the front side of a crystalline silicon solar cell, we formulate and analyse a two-dimensional mathematical model for electron flow across a poorly conducting (glass) layer situated between silver electrodes, based on the drift-diffusion (Poisson-Nernst-Planck) equations.  We devise and validate a novel spectral method to solve this model numerically.  We find that the current short-circuits through thin glass layer regions.  This enables us to determine asymptotic expressions for the average current density for two different canonical glass layer profiles.
\end{abstract}

% \begin{keywords}
% drift diffusion, asymptotic analysis, contact resistance, electrochemical systems, spectral methods, Poisson-Nernst-Planck
% \end{keywords}

% \begin{AMS}
% 30C20, 35B40, 35J47, 65N35, 78A35
% \end{AMS}

%\pagestyle{myheadings}
%\thispagestyle{plain}
%\markboth{Two-dimensional modelling of electron flow}{Two-dimensional modelling of electron flow}

\section{Introduction}

\begin{figure}
\begin{center}
 \scalebox{1}{
 \input{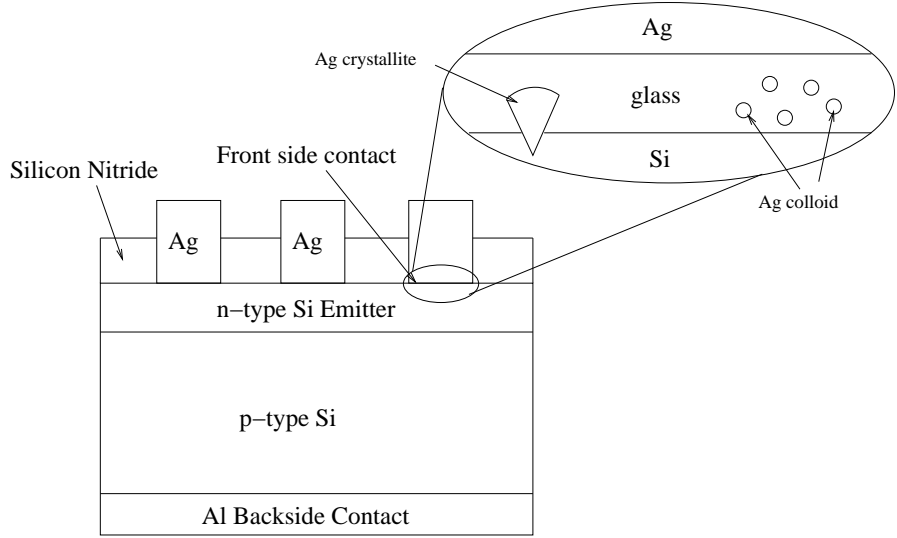}
}
\caption{Schematic diagram of a solar cell after firing \cite{Black13}. The detail shows the different possible geometries of silver inclusions: ``crystallite'' and ``colloid''.}
\label{fig:solarcell}
\end{center}
\end{figure}

Screen-printed crystalline silicon photovoltaic cells make up the majority of solar cells produced today. After manufacture \cite{Hong09}, the front contact of each cell consists of a silicon emitter and a silver electrode separated by a thin interfacial glass layer ($\sim10$ nm - $1\mu$m) that impedes electron flow \cite{Li09, Li11}. A schematic of the geometry is illustrated in figure \ref{fig:solarcell}.

The aim of this paper is to build upon the work outlined in \cite{Black13} to increase understanding of the local electron transport mechanisms at the front contact.  In the literature there is an ongoing debate over whether the presence of silver crystallites or of silver colloids in the glass layer is more effective in aiding electron transport \cite{Cabrera11, cheng2009nano,Ionkin11,Jeong10,Kontermann10,Kontermann11,Li11,Lin08,SchubertThesis}. In this paper we specifically focus on how the large variation in glass layer thickness affects the electron flow through the layer and consequently the effective contact resistance.  Since we are considering a poor conductor (glass) situated between two good conductors we expect the current to short circuit at thinner regions of the glass.  This feature will be investigated both numerically and through asymptotic calculations.  This study will also be of interest in adjacent fields since the situation in which a (relatively) poor conductor is 
positioned between two good 
conductors arises in many different physical applications,  for example in electrochemical thin-films \cite{Bazant05current, Bonnefort01, Chu2005electrochemical}, microelectronics \cite{cumberbatch2006nano, Ward90}, electrochromic glass \cite{mortimer1997electrochromic} and metal corrosion \cite{van10}.

We approach the problem in a similar manner to \cite{Black13}, where we analyse a mathematical model of electron flow based 
on the drift-diffusion (Poisson-Nernst-Planck) equations, omitting the 
boundary effects that 
are present at the silicon-glass and glass-silver interfaces.  The drift-diffusion equations are widely used to model current flow in a range of devices; see \cite{Please82,Markowich86,Markowich86b,Markowich84, murphy1992numerical, Richardson11, Bazant05current, Biesheuvel10, Chu2005electrochemical, Bonnefort01, Foster13, Ward90}.  These studies consider the drift-diffusion equations in one dimension; when considering the two-dimensional extension it is generally necessary to solve numerically. One of the first 2D schemes was devised by Slotboom \cite{Slotboom69} who reformulated the charge densities in terms of the quasi-fermi levels and solved using Gummel iteration \cite{gummel1964self} and the finite difference method.  In semiconductor devices, large gradients are often present near interfaces and therefore schemes using non-uniform meshes for finite difference methods have been formulated \cite{slotboom1973computer}.  For non-uniform meshes the finite-element method is more suitable and has been 
employed in a 
large number of 
studies, for example \cite{
adachi1979two} and more recently in \cite{Brinkman13,van10}.  A good review of the earlier numerical methods used is given by Snowden \cite{Snowden85}.  

Spectral methods are known to be good at resolving boundary layers and have been successfully employed in one dimension by Chu and Bazant \cite{Chu2005electrochemical} and through \textsf{Chebfun} \cite{chebfunv5, Tref2013} by Black et al. \cite{Black13} and Foster et al. \cite{Foster13}.  In this paper we formulate and employ a spectral method to solve the drift-diffusion equations in two dimensions.  The numerical method that we use is largely based upon the work of Birkisson and Driscoll \cite{Birkisson12}.  This involves using the Fr\'{e}chet derivatives of the problem, to implement Newton's method and thus obtain the numerical solution to a boundary value problem.    

In this paper we extend the model for electron transport through a glass layer formulated in \cite{Black13} to two dimensions.  After outlining two alternative mapping techniques to recast the governing equations in a rectangular domain, we describe our spectral numerical method.  By making use of asymptotic techniques, approximate expressions are found for the dependence of the resistance on the minimum thickness of the glass layer in two different canonical domain types.

\section{Mathematical model}

\subsection{Formulation of the problem}

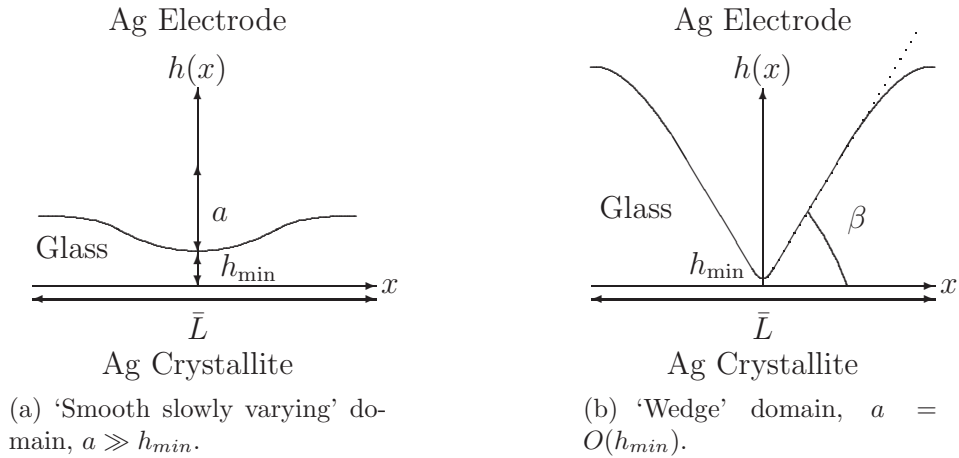
\begin{figure}%[htbp]
\centering
        \subfigure[`Smooth slowly varying' domain, $a \gg h_{min}$.]
               { \centering
  \begin{picture}(135,140)
%       % Drawing outline
	\put(5,35){\vector(1,0){130}}	% Bottom surface
	\qbezier(37,56.7)(67.5,39.7)(98,56.7) % Constructing top surface
	\qbezier(7.5,61.7)(30,61.7)(37.8,56.4)
	\qbezier(97.5,56.4)(105,61.7)(127.5,61.7)
	
	\put(67.5,35){\vector(0,1){12.5}} %hmin up vector
	\put(67.5,47.5){\vector(0,-1){12.5}} %hmin down vector
	\put(67.5,48.5){\vector(0,1){32.5}} % a up vector
	\put(67.5,80.5){\vector(0,-1){32.5}} % a down vector
	\put(67.5,50){\vector(0,1){60}} % h(x) vector
	\put(5,30){\vector(1,0){130}} % Lbar right vector
	\put(135,30){\vector(-1,0){130}} %Lbar left vector
	%\put(60,25){\vector(0,-1){15}}
	
	% Labelling Drawing
	\put(67.5,20){\makebox(0,0)[c]{$\bar{L}$}}
 	\put(67.5,135){\makebox(0,0)[c]{Ag Electrode}}
 	\put(67.5,5){\makebox(0,0)[c]{Ag Crystallite}}
 	\put(20,50){\makebox(0,0)[c]{Glass}}
 	\put(76,42){\makebox(0,0)[l]{$h_{\text{min}}$}}
 	\put(76,62.25){\makebox(0,0)[c]{$a$}}
 	\put(67.5,117){\makebox(0,0)[v]{$h(x)$}}
 	\put(137,35){\makebox(0,0)[l]{$x$}}
      \end{picture}
    \label{fig:smooth}
		} \qquad \qquad \qquad 
        \subfigure[`Wedge' domain, $a=O(h_{min})$.]
		{\centering
      \begin{picture}(130,140)
%       % Drawing outline
	\put(3,35){\vector(1,0){130}} % Bottom surface

	\put(72,41.6){\line(3,5){31}} % Wedge edge 1
	\put(64,41.6){\line(-3,5){31}} % Wedge edge 2
	\multiput(68,35)(2,3.3){30}{\line(0,1){0.3}} %Dashed line
	\qbezier(64,41.6)(68,34)(72,41.6) % Curving the wedge

	\qbezier(33,93)(15,120)(3,118) % Curving off the left side
	\qbezier(103,93)(121,120)(133,118) % Curving off the right side
	\arc(100,35){17.2} % Labelling angle
	\put(3,30){\vector(1,0){130}} % Lbar left vector
	\put(133,30){\vector(-1,0){130}} %Lbar right vector
	\put(68,35){\vector(0,1){75}} % h(x) vector
	
	% Labelling Drawing
	\put(68,20){\makebox(0,0)[c]{$\bar{L}$}}
	\put(68,135){\makebox(0,0)[c]{Ag Electrode}}
 	\put(68,5){\makebox(0,0)[c]{Ag Crystallite}}
 	\put(20,65){\makebox(0,0)[c]{Glass}}
 	\put(40,43){\makebox(0,0)[l]{$h_{\text{min}}$}}
 	\put(100,60){\makebox(0,0)[l]{$\beta$}}
 	\put(68,117){\makebox(0,0)[c]{$h(x)$}}
 	\put(135,35){\makebox(0,0)[l]{$x$}}

      \end{picture}
     \label{fig:wedge}
		}
  \caption{The two different domain shapes that are investigated; $h_{min} =$ minimum glass thickness, $a=$ radius of curvature and $\beta=$ angle from horizontal.}
    \label{fig:smoothandwedge}
\end{figure}

We consider two-dimensional conduction through a symmetric periodic glass layer of wave-length $\bar{L}$ that has varying thickness, $h(x)$, driven by an applied electric potential difference $\bar{\Phi}$.  The glass layer has average thickness $H$ and minimum thickness $h_{\text{min}}$. We focus on a system where a glass layer separates a silver crystallite from the silver electrode.  We want to investigate the short circuiting of the current through thinner regions of the glass layer.  We therefore consider two different generic cases for the variation of the thickness of the glass layer: (i) the radius of curvature, $a$, of the silver electrode is much larger than the minimum thickness of the glass layer $h_{\text{min}}$, which implies the thickness of the glass layer varies slowly, see figure~\ref{fig:smooth}; (ii) the radius of curvature $a = O(h_{\text{min}})$ so that the domain appears wedge-like and the thickness of the glass layer varies quickly, see figure \ref{fig:wedge}. In the latter 
case we find it more convenient to parametrise the domain by the angle 
from the horizontal, $\beta$.

\subsection{Nondimensional model}
\label{sec:nondim}
\begin{figure}
 \centering
 \scalebox{1}{
\begin{picture}(240,225)
 % Drawing outline
  \put(0,45){\vector(1,0){240}}
  \put(15,30){\vector(0,1){180}}
  \put(225,30){\line(0,1){180}}
  \qbezier(75,150)(120,15)(165,150)
  \qbezier(15,180)(60,183)(75,150)
  \qbezier(165,150)(180,183)(225,180)
  
 % Labelling yaxis
 \put(0,45){\makebox(0,0)[r]{$0$}}
 \put(0,135){\makebox(0,0)[r]{$\phi_x = 0$}}
 \put(0,105){\makebox(0,0)[r]{$n_x = 0$}}
 %\put(-12,80){\makebox(0,0)[r]{$H$}}
 \put(15,225){\makebox(0,0)[c]{$y$}}

 %Labelling top surface
 \put(120,175){\makebox(0,0)[c]{Silver Electrode}}
 \put(195,155){\makebox(0,0)[c]{$y=h(x)$}}
 \put(120,155){\makebox(0,0)[c]{$n=1$}}
 \put(120,125){\makebox(0,0)[c]{$\phi=0$}}

 %Labelling Glass
 \put(121,63.75){\makebox(0,0)[l]{\small{$h_{\text{min}}$}}}
 \put(120,45){\vector(0,1){37.5}}
 \put(120,82.5){\vector(0,-1){37.5}}
 \put(60,80){\makebox(0,0)[c]{Glass}}
 
 %Labelling RHS
 \put(240,135){\makebox(0,0)[l]{$\phi_x = 0$}}
 \put(240,105){\makebox(0,0)[l]{$n_x = 0$}}
 
 %Labelling x axis
 \put(15,22.5){\makebox(0,0)[c]{$-L/2$}}
 \put(90,30){\makebox(0,0)[c]{$\phi=\Phi$}}
 \put(150,30){\makebox(0,0)[c]{$n=1$}}
 \put(225,22.5){\makebox(0,0)[c]{$L/2$}}
 \put(255,45){\makebox(0,0)[c]{$x$}}
 \put(120,7){\makebox(0,0)[c]{Silver Crystallite}}

 \end{picture}
}
\caption{Schematic of physical domain and boundary conditions.}
\label{fig:schemphysical}
 \end{figure}
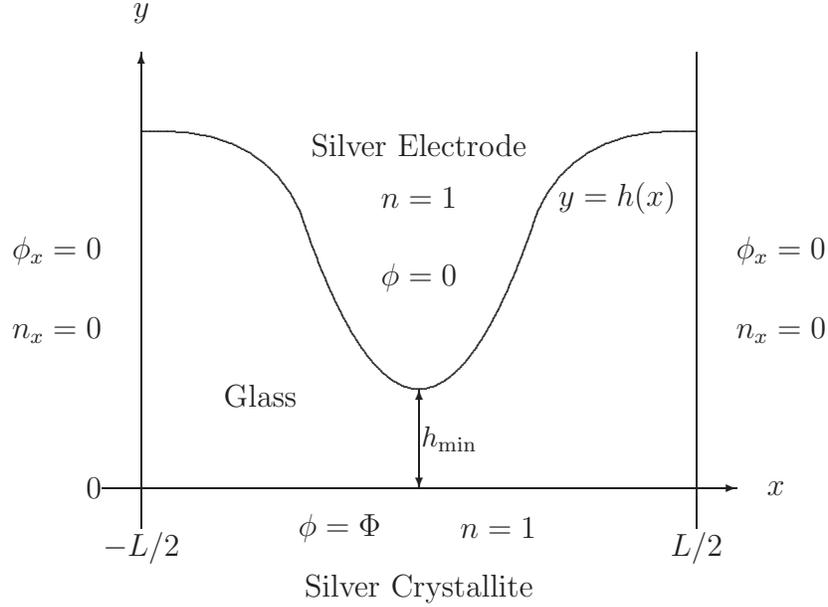
 To formulate the mathematical model we make the following assumptions; see \cite{Black13} for further details.
\begin{enumerate}
 \item Charge is predominantly carried by electrons and not holes.
 \item The system is in quasi-steady state. 
 \item The flux $\ve{j}$ of electrons results from combination of drift and diffusion.
 \item The electric potential $\phi$ satisfies Poisson's equation.
\end{enumerate}
We then follow the same nondimensionalisation as in \cite{Black13} to obtain the following governing equations
\begin{align}
 \nabla \cdot \ve{j} &= 0, \label{eq:consvj}\\
  \ve{j} &= -n \nabla \phi -\nabla n, \label{eq:jdrifdif}  \\
 \nu^2\nabla^2 \phi &= -n, \label{eq:jpois} 
\end{align}
where $n$ is the density of free electrons, and the boundary
conditions are
\begin{align}
  \phi(x,h(x)) & = 0,&  n(x, h(x)) & = 1, & \phi(x,0) & = \Phi, & n(x,0) & = 1,\nonumber\\
 \phi_x(L/2,y) & = 0, & n_x(L/2,y) & = 0, & \phi_x(-L/2,y) & = 0, & n_x(-L/2,y) & = 0,  \label{eq:bound}
\end{align}
where we assume periodic boundary conditions in the $x$-direction. A schematic of the domain and the boundary conditions is given in figure \ref{fig:schemphysical}.  The dimensionless parameters are
\begin{align}
 L & = \f{\bar{L}}{H}, & 
\nu &= \frac{H_D}{H}=\sqrt{\frac{\ep k_B T}{q^2 H^2 n_1}}, & \Phi & =  \f{q \bar{\Phi}}{k_B T},
\end{align}
where $L$ is the nondimensional length of the glass layer and can be thought of as the aspect ratio ($L:1$), $H$ is the average thickness of the glass layer,  $H_D$ is the Debye length, $\epsilon$ is the absolute permittivity, $k_B$ is Boltzmann's constant, $T$ is absolute temperature, $n_1$ is the electron density at the silver electrode interface and $q$ is the charge on an electron. We expect the Debye length typically to be small compared with the average glass layer thickness, and hence the parameter $\nu$ to be small. It is when the Debye length is small compared with the average glass layer thickness that the short-circuting effects investigated are especially evident. We only consider cases where $\Phi>0$ since this corresponds to electron flow into the silver electrode.

\subsection{Model outputs}
From the solution of the mathematical model there are a number of different quantities we will calculate to gain greater insight into the electron transport across the glass layer.  Firstly, we will determine the average current density through the glass layer which is given by
\begin{equation}
 Q = \f{1}{L}\int_{-L/2}^{L/2}\left[-n\pd{\phi}{y} - \pd{n}{y}\right]_{y=0}\, \sd x; \label{eq:Q}
\end{equation}
  the total current is given by $QL$.  To investigate qualitatively the short circuiting of the current through thinner glass layer regions we calculate: (i) the electron trajectories through the glass layer; (ii) the normalised cumulative current
\begin{equation}
 \Cc(x) = \f{2}{QL}\int_{0}^x \left[-n\pd{\phi}{y} - \pd{n}{y}\right]_{y=0} \,\sd x. \label{eq:cumcurr}
\end{equation}
The deviation of the function $\Cc(x)$ from linear indicates focusing of the current.

Finally, we define the effective resistance of the glass layer as 
\begin{equation}
 R = \f{\Phi}{Q}.
\end{equation}

\section{Transformation to rectangular domain}
\label{sec:Transformation}

\subsection{Introduction}

In this paper we solve (\ref{eq:consvj})--(\ref{eq:bound}) using a spectral numerical method.  The numerical solution is facilitated by mapping the domain to a rectangle.    In this section we outline two alternative mapping techniques: a hodograph formulation, particularly suitable for investigating the slowly varying domain (figure \ref{fig:smooth}) and a conformal map formulation \cite{bazant2004conformal}, ideal for investigating the wedge domain (figure \ref{fig:wedge}).

\subsection{Hodograph transformation}

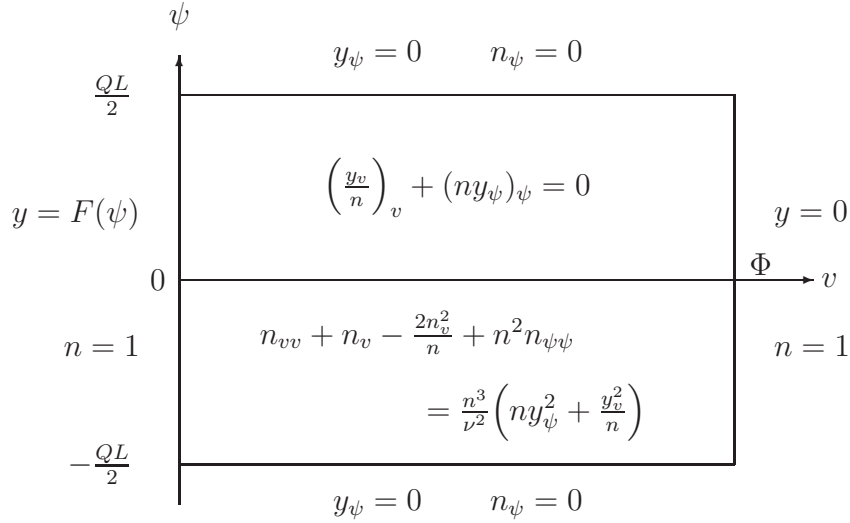
\begin{figure}
\centering
\scalebox{1}{
      \begin{picture}(240,200)
	% Drawing outline
	  \put(18,92.5){\vector(1,0){240}}
	  \put(18,7.5){\vector(0,1){170}}
	  \put(228,22.5){\line(0,1){140}}
	  \put(18,162.5){\line(1,0){210}}
	  \put(18,22.5){\line(1,0){210}}
	  
	% Labelling \psi axis
	\put(0,22.5){\makebox(0,0)[r]{$-\f{QL}{2}$}}
	\put(3,117.5){\makebox(0,0)[r]{$y = F(\psi)$}}
	\put(3,67.5){\makebox(0,0)[r]{$n=1$}}
	\put(0,162.5){\makebox(0,0)[r]{$\f{QL}{2}$}}
	\put(18,192.5){\makebox(0,0)[c]{$\psi$}}

	%Labelling top surface
	\put(93,177.5){\makebox(0,0)[c]{$y_{\psi} = 0$}}
	\put(153,177.5){\makebox(0,0)[c]{$n_{\psi}=0$}}
	
	%Labelling RHS
	\put(243,117.5){\makebox(0,0)[l]{$y = 0$}}
	\put(243,67.5){\makebox(0,0)[l]{$n = 1$}}
	
	%Labelling v axis
	\put(10,92.5){\makebox(0,0)[c]{$0$}}
	\put(93,7.5){\makebox(0,0)[c]{$y_{\psi} = 0$}}
	\put(153,7.5){\makebox(0,0)[c]{$n_{\psi} = 0$}}
	\put(238,94.5){\makebox(0,0)[b]{$\Phi$}}
	\put(264,92.5){\makebox(0,0)[c]{$v$}}
	
	% Adding in Equations
	\put(123,127.5){\makebox(0,0)[c]{$\Bigl(\f{y_v}{n}\Bigr)_v + (ny_{\psi})_{\psi} = 0$}}
	\put(108,72.5){\makebox(0,0)[c]{$n_{vv} + n_v - \f{2n_v^2}{n} + n^2n_{\psi\psi}$}}
	\put(153,42.5){\makebox(0,0)[c]{$= \f{n^3}{\nu^2}\Bigl(ny_{\psi}^2 + \f{y_v^2}{n}\Bigr)$}}
      \end{picture}
      }
      \caption{Schematic of hodograph plane domain and boundary conditions.}
      \label{fig:schem2b}
\end{figure}

Equation (\ref{eq:consvj}) implies the existence of a streamfunction $\psi(x,y)$ such that
\begin{equation}
 \ve{j} = (\psi_y, -\psi_x).
\end{equation}
The introduction of the electro-chemical potential $v = \phi + \log n$ enables us to write (\ref{eq:jdrifdif}) as 
\begin{align}
 \psi_y & = -nv_x, & \psi_x & = nv_y. \label{eq:Cr}
\end{align}
Now by considering the boundary conditions (\ref{eq:bound}), we observe that the change of independent variables from $(x,y) \rightarrow (v,\psi)$, maps the solution domain onto the rectangle $(v,\psi) \in [0,\Phi] \times [-QL/2,QL/2]$.  The Jacobian of the transformation is given by 
 \begin{equation}
  J = \pd{(v,\psi)}{(x,y)} = v_x\psi_y - v_y\psi_x.\label{eq:Jacobian}
 \end{equation}
Now by inverting the relations (\ref{eq:Cr}) we find that
\begin{align}
 x_v & = -ny_{\psi}, & y_v & = nx_{\psi} \label{eq:Crinvert}
\end{align}
and hence
\begin{equation}
 \f{1}{J}  = x_vy_{\psi} - x_{\psi}y_v  =  - n\left(ny_{\psi}^2 + \f{y_v^2}{n}\right). \label{eq:Jinv}
\end{equation}
To formulate our problem in the Hodograph plane we want to write (\ref{eq:jpois}) in terms of $v$ and $\psi$. We use the chain rule to find that
 \begin{align}
   \psi_{xx} + \psi_{yy} & = Jn\pd{\phi}{\psi}, \label{eq:psixx} \\
   v_{xx} + v_{yy} & = \f{J}{n} \left(1 - \pd{\phi}{v}\right) \label{eq:vxx}.
 \end{align}
 Substituting (\ref{eq:psixx}) and (\ref{eq:vxx}) into (\ref{eq:jpois}), and making use of (\ref{eq:Cr}) and (\ref{eq:Jacobian}) we obtain
 \begin{equation}
 \left(\f{\phi_v}{n}\right)_v + \left(n\phi_{\psi}\right)_{\psi} = \f{n}{\nu^2J}. \label{eq:phinhod}
\end{equation}
 Finally, substituting (\ref{eq:Jinv}) into (\ref{eq:phinhod}) and eliminating $x$ from (\ref{eq:Crinvert}) we find the governing equations in the hodograph plane:
\begin{align}
 n_{vv} + n_v - \f{2n_v^2}{n} + n^2n_{\psi\psi} & = \f{n^3}{\nu^2}\left(n y_{\psi}^2 + \f{y_v^2}{n}\right), \\
 \left(\f{y_v}{n}\right)_v + (ny_{\psi})_{\psi} & = 0,
\end{align}
with boundary conditions
\begin{align}
 n(0,\psi) & = 1, & n(\Phi,\psi) & = 1, & y(0,\psi) &=F(\psi), & y(\Phi,\psi) & =0 \\
  n_{\psi}(v,-QL/2)& = 0, & n_{\psi}(v,QL/2) & = 0, & y_{\psi}(v,-QL/2) & = 0, & y_{\psi}(v,QL/2) & = 0.
\end{align}
A schematic of the hodograph plane and boundary conditions is shown in figure \ref{fig:schem2b}.  We note that in this formulation the method has become semi-inverse in nature as we now specify the total current $QL$ and solve for the length of the domain~$L$. We also specify the function $y(0,\psi) = F(\psi)$ that determines indirectly the shape of the top surface in the physical plane, $y = h(x)$. 

To plot our results in the physical plane we  calculate $x(v,\psi)$:
\begin{equation} \label{eq:xhod}
 x(v,\psi) = -\int_{-QL/2}^{\psi} \f{y_v(v,\psi')}{n(v,\psi')} \quad \sd \psi'.
\end{equation}
 In the hodograph plane the length of the physical domain is now an output and given by
\begin{equation}\label{eq:Lhod}
 L(v) =  -\int_{-QL/2}^{QL/2}\f{y_v(v,\psi)}{n(v,\psi)} \quad \sd \psi;
\end{equation}
$L$ should be independent of $v$ and any variation in $L$ over the range of $v$ is due to numerical error.  The electron trajectories are easily found as curves of constant $\psi$.  Finally, the cumulative current is found by inverting the function $x(\Phi,\psi)$ to determine $\psi(x,0)$.

\subsection{Conformal mapping}

\subsubsection{Transformed problem}

\begin{figure}
 \centering 
 \scalebox{1}{
      \begin{picture}(360,170) 
      % Drawing outline 
	\put(3,22.5){\vector(1,0){350}}
	\put(18,22.5){\line(0,1){120}}
	\put(178,22.5){\line(0,1){27}}
	\put(178,65){\line(0,1){32}}
	\put(178,115){\vector(0,1){47}}
	%\multiput(178,22.5)(0,10){13}{\line(0,1){5}}
	
	%\put(178,22.5){\vector(0,1){140}}
	\put(338,22.5){\line(0,1){120}}
	\put(18,142.5){\line(1,0){320}}
	
      % Labelling \eta axis 
      \put(178,11){\makebox(0,0)[c]{$0$}}
      \put(3,97.5){\makebox(0,0)[r]{$\phi_{\xi} = 0$}}
      \put(3,47.5){\makebox(0,0)[r]{$n_{\xi} = 0$}}
      \put(178,135){\makebox(0,0)[r]{$\eta^*$}}
      \put(178,166.5){\makebox(0,0)[b]{$\eta$}}

      %Labelling top surface 
      \put(90.5,152.5){\makebox(0,0)[c]{$\phi = 0$}}
      \put(265.5,152.5){\makebox(0,0)[c]{$n=1$}}
      
      %Labelling RHS 
      \put(348,97.5){\makebox(0,0)[l]{$\phi_{\xi} = 0$}}
      \put(348,47.5){\makebox(0,0)[l]{$n_{\xi} = 0$}}
      
      %Labelling \xiaxis 
      \put(18,11){\makebox(0,0)[c]{$-L/2$}}
      \put(90.5,11){\makebox(0,0)[c]{$\phi=\Phi$}}
      \put(265.5,11){\makebox(0,0)[c]{$n=1$}}
      \put(338,11){\makebox(0,0)[c]{$L/2$}}
      \put(358,22.5){\makebox(0,0)[c]{$\xi$}}
      
      % Adding in Equations 
      \put(178,57.5){\makebox(0,0)[c]{$\nu^2 \left(\phi_{\xi\xi}
+\phi_{\eta\eta} \right)  = -G(\xi,\eta)n$}}
      \put(168,107.5){\makebox(0,0)[c]{$\nu^2\left(
n_{\xi\xi} + n_{\eta\eta} + n_\xi\phi_\xi+n_\eta\phi_\eta \right)$}}
      \put(212,87.5){\makebox(0,0)[c]{$=G(\xi,\eta)n^2$}}
      \end{picture}
      }
    \caption{Schematic of conformally mapped domain and boundary
      conditions.}
    \label{fig:schem2a}
\end{figure}
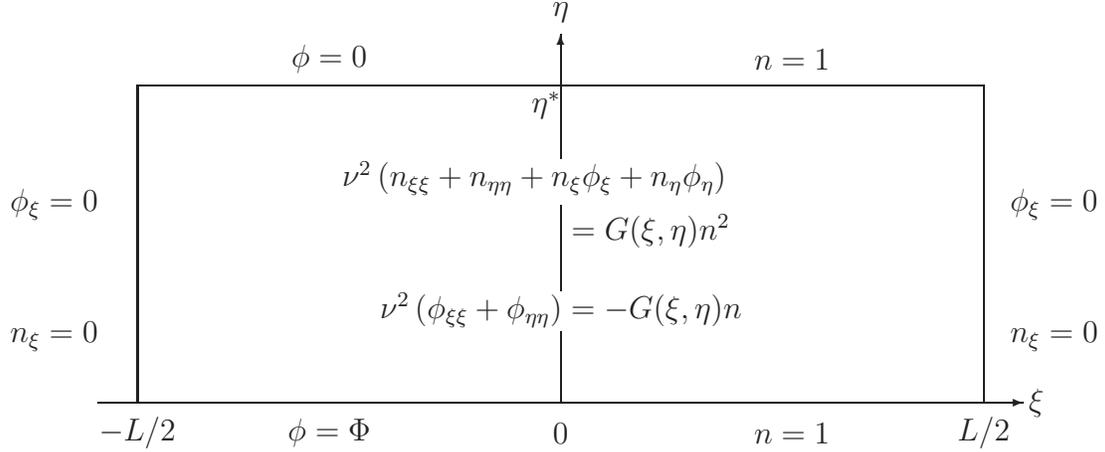

As an alternative to the hodograph transformation described above,
suppose instead that the rectangle $R=[-L/2,L/2]\times[0,\eta^*]$ is
mapped onto the physical domain of interest
$\{(x,y):0\leq y\leq h(x);-L/2\leq x\leq L/2\}$ by a conformal
transformation
\begin{equation}\label{zetatoz}
x+\ii y=z=f(\xi+\ii\eta)=f(\zeta).
\end{equation}
Application of (\ref{zetatoz}) to (\ref{eq:consvj})--(\ref{eq:jpois})
gives the following equations in the transformed frame:
\begin{align}
 \nu^2 \left(\pdd{\phi}{\xi} + \pdd{\phi}{\eta}\right) + G(\xi,\eta) n& = 0, \label{eq:2Dpoisconf}\\
  \pdd{n}{\xi} + \pdd{n}{\eta} + \pd{n}{\xi}\pd{\phi}{\xi} + \pd{n}{\eta}\pd{\phi}{\eta} - G(\xi,\eta)\f{n^2}{\nu^2} & = 0, \label{eq:2Ddiffconf}
\end{align}
where $G(\xi,\eta)=|f'(\xi+\ii\eta)|^2$, and the boundary conditions
(\ref{eq:bound}) become
 \begin{align}
  \phi(\xi,0) & = \Phi, & n(\xi,0) & = 1, & \phi(\xi,\eta^*) &=0, & n(\xi,\eta^*) & =1 \nonumber\\
  \phi_{\xi}\left(-\f{L}{2},\eta\right)& = 0, & n_{\xi}\left(-\f{L}{2},\eta\right) & = 0, & \phi_{\xi}\left(\f{L}{2},\eta\right) & = 0, & n_{\xi}\left(\f{L}{2},\eta\right) & = 0. \label{eq:2Dconfbound}
 \end{align}
A schematic of the domain and boundary conditions is shown in
figure~\ref{fig:schem2a}.
The height $\eta^*$ of the rectangular domain $R$ is chosen to
ensure that the nondimensional average thickness of the glass layer is 
equal to $1$, that is,
\begin{equation}
 \text{Area} = - \f{1}{2} \im \left[\int_{-L/2}^{L/2} f(\xi -\ii\eta^*) 
     f'(\xi +\ii\eta^*) \, \sd \xi \right]+ \f{L}{2}
y_{\text{max}}\equiv L,
  \end{equation}
where $y_{\text{max}}$ is the maximum value that $y$ takes in the $z$
plane.

The calculation of the average current density $Q$ is achieved using
the relation
\begin{equation}
  Q = \f{1}{L}\int_{-L/2}^{L/2} j^{\eta}(\xi,0) \, \sd \xi,
\label{eq:Qconf}
\end{equation}
where $j^{\eta}(\xi,\eta) = - n\phi_{\eta} - n_{\eta}$. To obtain the
electron trajectories, we solve the equations
\begin{align}
  \od{\xi}{t} & = - \left(n \pd{\phi}{\xi} + \pd{n}{\xi}\right),
&
  \od{\eta}{t} & = -\left(n\pd{\phi}{\eta} + \pd{n}{\eta}\right),
 \end{align}
using the overloaded \textsf{ode45} function in
\textsf{Chebfun} and then map back to the physical
domain. Finally, the cumulative current plots can be produced by
noting the relation
\begin{equation}
 \Cc(x) =\f{2}{QL}\int_{0}^{x} j^{\eta}(\xi,0) \pd{\xi}{x} \sd x,
\end{equation}
and using the overloaded \textsf{cumsum} command in
\textsf{Chebfun}. 

\subsubsection{Mapping function}

\begin{figure}[hbtp!]
  \subfigure{
    \begin{picture}(210, 100)(-18,0)     
      \put(-5,24){\makebox(0,0)[c]{\tiny{$0$}}}
      \put(-10,80){\makebox(0,0)[c]{\tiny{$\eta_{max}$}}}
      \put(-5,52){\makebox(0,0)[c]{\tiny{$\eta^*$}}}
      \includegraphics[width=0.432\textwidth]{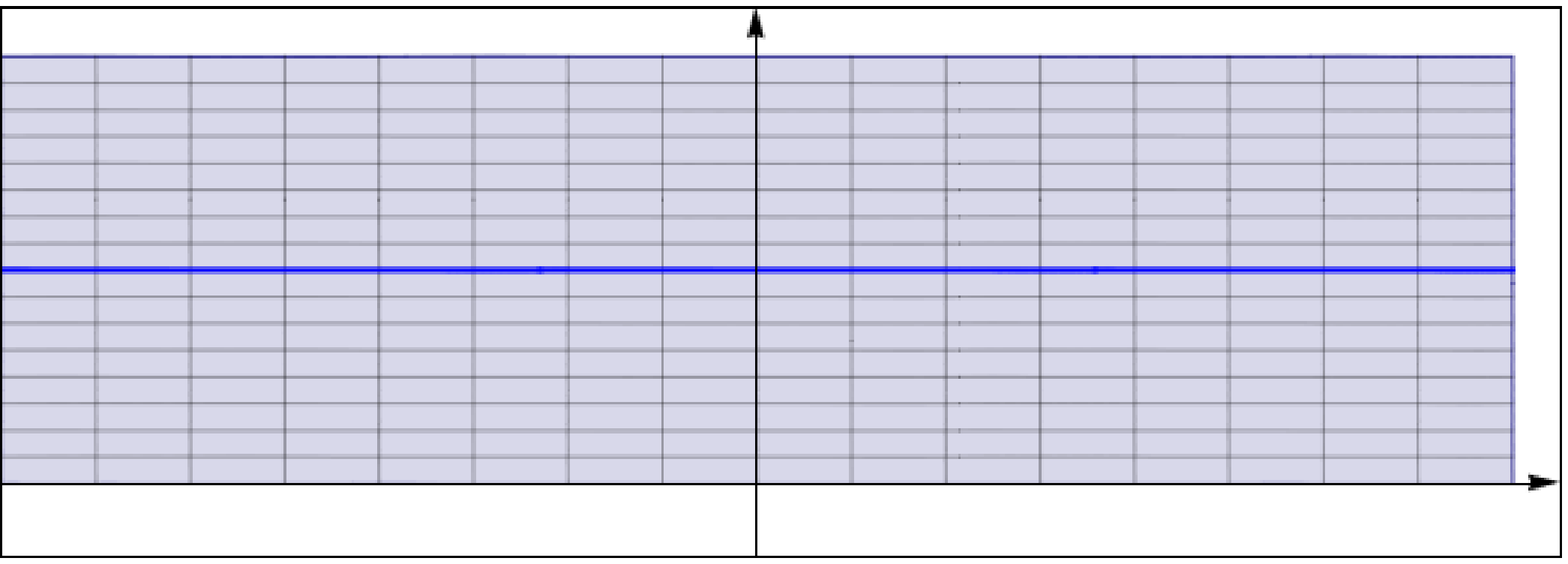}
      \put(-201,8){\makebox(0,0)[c]{\tiny{$-\pi$}}}
      \put(-105,8){\makebox(0,0)[c]{\tiny{$0$}}}
      \put(-5,8){\makebox(0,0)[c]{\tiny{$\pi$}}}
      
      \put(-105,90){\makebox(0,0)[c]{\tiny{$\eta$}}}
      \put(4,24){\makebox(0,0)[c]{\tiny{$\xi$}}}

    \end{picture}
  }
  \qquad
  \subfigure{
    \begin{picture}(210,100)(0, 0)
      \put(-5,24){\makebox(0,0)[c]{\tiny{$0$}}}
      \put(-5,40){\makebox(0,0)[c]{\tiny{$\varepsilon$}}}
      \includegraphics[width=0.432\textwidth]{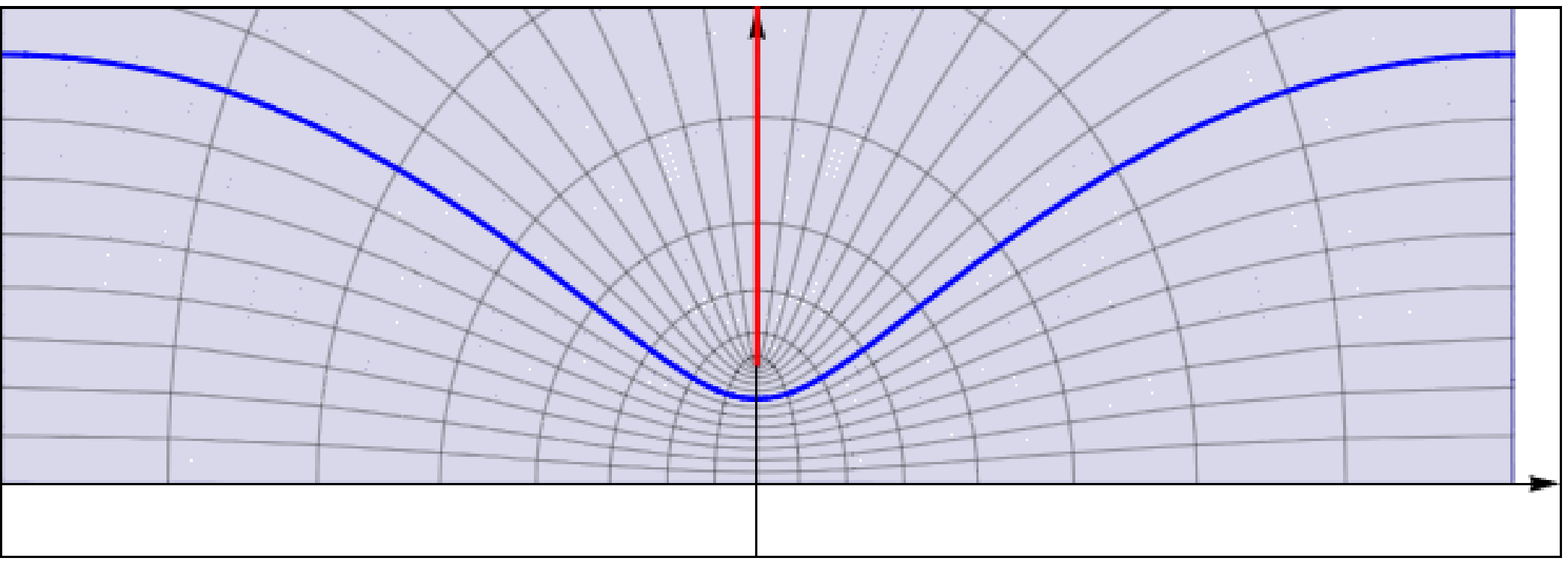}
      \put(-201,8){\makebox(0,0)[c]{\tiny{$-\pi$}}}
      \put(-105,8){\makebox(0,0)[c]{\tiny{$0$}}}
      \put(-5,8){\makebox(0,0)[c]{\tiny{$\pi$}}}
      
      \put(-105,90){\makebox(0,0)[c]{\tiny{$y$}}}
      \put(4,24){\makebox(0,0)[c]{\tiny{$x$}}}
    \end{picture}
  }
  \caption{The conformal map (\ref{eq:confmap}) with $L=2\pi$ and $\varepsilon=0.5$, maps
    the $\zeta$-plane on the left to the $z$ plane on the right.}
  \label{fig:Confmapeg}
\end{figure}
For a given physical domain a numerical conformal map can be constructed by making use of the Schwarz-Christoffel transformation.  However, to explore domains with a wedge-like geometry, we use the following mapping function, see Hale \& Tee
\cite{hale2009conformal}:
 \begin{equation}\label{eq:confmap}
f(\zeta) =  \f{L}{\pi}\arcsin\left(
\frac{\displaystyle\tanh\left({\epsilon}/{2}\right)
\sn\left({2\K\zeta}/{L}\right)}
{\displaystyle\dn\left({2\K\zeta}/{L}\right)}\right), 
\end{equation}
where $\K$ denotes the complete elliptic integral of the first kind;
$\sn$ and $\dn$ are Jacobi elliptic functions. The suppressed argument
$m$ of $\K(m)$, $\sn(u|m)$ and $\dn(u|m)$ is defined by
\begin{equation}
m=\sech^2\left({\epsilon}/{2}\right).
\end{equation}
Hence the map (\ref{eq:confmap}) depends on one parameter $\epsilon$
in addition to the length $L$ of the periodic domain.  

The function $f(\zeta)$ is univalent on the rectangle
$R=[-L/2,L/2]\times[0,\eta^*]$ provided
\mbox{$\eta^*<\eta_{\text{max}}$}, where
\begin{equation}
\eta_{\text{max}}=\frac{L}{2}\,\frac{\K(1-m)}{\K(m)}.
\end{equation}
As shown in figure~\ref{fig:Confmapeg}, the function
(\ref{eq:confmap}) maps the rectangle
$[-L/2,L/2]\times[0,\eta_{\text{max}}]$
in the $\zeta$-plane to a strip $[-L/2,L/2]\times[0,\infty]$ in the
$z$-plane, minus a branch cut from $\ii\epsilon L/2\pi$ to
$\ii\infty$. For $\eta^*<\eta_{\text{max}}$, the line
$\im(\zeta)=\eta^*$ is mapped to a periodic curved upper surface that
is wrapped around the branch cut. Moreover, when $\epsilon$ and
$\eta^*$ are small, the geometry approaches that of a rounded wedge
which narrowly avoids intersecting the lower surface $y=0$ at $x=0$.
In the limit as $\epsilon\rightarrow0$, we find that
\begin{align}\label{eq:confmapinn}
\eta_{\text{max}}&\sim\frac{\pi L}{4\log(8/\epsilon)},
&
z=f(\zeta)&\sim\frac{\epsilon L}{2\pi}\sinh\left(
\frac{2\zeta\log(8/\epsilon)}{L}\right).
\end{align}
The parameters $\eta^*$ and $\epsilon$ may thus be related to the
asymptotic minimum layer thickness $h_{\text{min}}$ and wedge angle
$\beta$ using the formulae
\begin{align}
\eta^*&=\frac{\beta L}{2\log(8/\epsilon)},
&
\epsilon&=\frac{2\pi h_{\text{min}}}{L\sin\beta},
\end{align}
and the local geometry is that of a hyperbola, with
\begin{equation}\label{hyperbola}
y^2\cos^2\beta-x^2\sin^2\beta=h_{\text{min}}^2\cos^2\beta.
\end{equation}
The use of the analytic conformal map, (\ref{eq:confmap}), rather than a numerical Schwarz-Christoffel map, reduces computational effort and allows easy parameterisation of a range of wedge-like domains.

\subsection{Summary}
\label{sec:transsum}

We have formulated two different techniques for mapping the governing equations to a rectangular domain.  These approaches each have advantages and disadvantages and are appropriate for treating different domains. The conformal map formulation is straightforward to apply and, if an analytic map can be found as shown above, allows for the easy manipulation of the domain shape.  Of course in general the Schwarz-Christoffel toolbox can be used to obtain numerical mapping functions, although this adds computational effort \cite{driscoll2005algorithm}.  The calculation of the electron trajectories and cumulative current also requires significant additional numerical calculations.  In contrast, the hodograph plane formulation makes it easy to determine these desired quantities.  In theory, the hodograph plane formulation allows the treatment of a wide range of domain shapes. However, unless \emph{a priori} knowledge of an appropriate $y=F(\psi)$ is known, it is difficult to reproduce a given domain.

\section{Numerical method}

\subsection{Introduction}

Via either of the transformations described in Section \ref{sec:Transformation}, we may henceforth assume that the governing equations are posed on a rectangular domain making them more suitable for numerical solution.  The numerical method we will use is a 2D spectral method.  The basis of spectral methods is to take discrete data on a grid, interpolate this data with a global function and then evaluate the derivative of the interpolating function on the grid. Given a periodic problem one would typically use trigonometric interpolants on equispaced points, and given a non-periodic problem it is commonplace to use Chebyshev polynomial interpolants on Chebyshev spaced points \cite{trefethen2000spectral}. Both formulations laid out in Section \ref{sec:Transformation} are periodic in one direction and non-periodic in the other, and we will therefore choose the appropriate interpolants in each direction accordingly.

\subsection{Demonstration of implementation}

We will give a detailed description of how the numerical method is implemented on the conformally mapped problem,  (\ref{eq:2Dpoisconf})--(\ref{eq:2Dconfbound}); this numerical method can be applied to analagous nonlinear elliptic equations in a rectangle and therefore exactly the same methodology is used for the hodograph plane formulation.  To apply Newton's method we calculate the Fr\'{e}chet derivatives of the governing equations with respect to the two unknowns $\phi$ and $n$:
\begin{align}
 \pd{(\ref{eq:2Dpoisconf})}{\phi}: & \quad \nu^2 \left(\pdd{}{\xi} + \pdd{}{\eta}\right), \\
 \pd{(\ref{eq:2Dpoisconf})}{n}: & \quad G(\xi,\eta), \\
 \pd{(\ref{eq:2Ddiffconf})}{\phi}: & \quad \pd{n}{\xi}\pd{}{\xi} + \pd{n}{\eta}\pd{}{\eta},\\
 \pd{(\ref{eq:2Ddiffconf})}{n}: & \quad \pdd{}{\xi} + \pdd{}{\eta} + \pd{\phi}{\xi}\pd{}{\xi} + \pd{\phi}{\eta}\pd{}{\eta} - \f{2n}{\nu^2}G(\xi,\eta).
\end{align}
Given an approximation $(\phi_k, n_k)$ to the solution, we therefore calculate an improved approximation $(\phi_{k+1}, n_{k+1}) = (\phi_k, n_k) + \gamma_k(u_k^{\phi}, u_k^n)$, where $u_k^{\phi}$, $u_k^n$ are the updates and $\gamma_k$ is the damping parameter that increases the chance of an initial guess converging to a solution \cite{Birkisson12}.  The updates satisfy the linear partial differential equations
\begin{align}
  \begin{pmatrix}
   \nu^2 \left(\pdd{}{\xi} + \pdd{}{\eta}\right) & \quad G(\xi,\eta) \\
   \pd{n_k}{\xi}\pd{}{\xi} + \pd{n_k}{\eta}\pd{}{\eta} & \pdd{}{\xi} + \pdd{}{\eta} + \pd{\phi_k}{\xi}\pd{}{\xi} + \pd{\phi_k}{\eta}\pd{}{\eta} - \f{2n_k}{\nu^2}G(\xi,\eta)
  \end{pmatrix}
  \begin{pmatrix}
   u_k^{\phi} \\
   u_k^n
  \end{pmatrix} \nonumber \\
  =-
  \begin{pmatrix}
   \nu^2 \left(\pdd{\phi_k}{\xi} + \pdd{\phi_k}{\eta}\right) + G(\xi,\eta) n_k, \\
   \pdd{n_k}{\xi} + \pdd{n_k}{\eta} + \pd{n_k}{\xi}\pd{\phi_k}{\xi} + \pd{n_k}{\eta}\pd{\phi_k}{\eta} - G(\xi,\eta)\f{n_k^2}{\nu^2}
  \end{pmatrix}. \label{eq:symsystem}
\end{align}
We ensure that the initial guess $(\phi_0, n_0)$ satisfies the boundary conditions in the $\eta$ direction so we can specify homogenous boundary conditions for the update $u$. We also specify periodic boundary conditions in the $\xi$ direction, i.e.
\begin{align}
 u^{\phi}_{\xi}(0,\eta) & = 0, & u^{\phi}_{\xi}(L,\eta) & = 0, &
 u^{\phi}(\xi,0) & = 0, & u^{\phi}(\xi,\eta^*) & = 0, \notag \\
 u^{n}_{\xi}(0,\eta) & = 0, & u^{n}_{\xi}(L,\eta) & = 0, &
 u^{n}(\xi,0) & = 0, & u^{n}(\xi,\eta^*) & = 0. \label{eq:numbound}
\end{align}
The successive approximations are calculated repeatedly until the update $||u_k|| < 10^{-5}$.  

\begin{figure}[hbtp!]
  \begin{minipage}{0.46\textwidth}
    \begin{equation}
    \phi_k = \begin{pmatrix}
	      \phi_k^1 \\
	      \phi_k^2 \\
	      \phi_k^3 \\
	      \vdots \\
	      \phi_k^{N_{\eta}N_{\xi}}
	      \end{pmatrix} \quad \text{where} \nonumber
    \end{equation}

    \end{minipage}
    \begin{minipage}{0.46\textwidth}
    \scalebox{0.9}{
      \begin{picture}(180,180)
    %       % Drawing outline
	    \put(0,10){\vector(1,0){165}}
	    \put(0,10){\vector(0,1){165}}
	    \put(157.5,10){\line(0,1){157.5}}
	    \put(0,167.5){\line(1,0){157.5}}
	    
	    % Putting on nodes
	    \put(0,10){\circle{4}}
	    \put(0,33){\circle{4}}
	    \put(0,88.75){\circle{4}}
	    \put(0,144.4){\circle{4}}
	    \put(0,167.5){\circle{4}}
	    
	    \put(39.375,10){\circle{4}}
	    \put(39.375,33){\circle{4}}
	    \put(39.375,88.75){\circle{4}}
	    \put(39.375,144.4){\circle{4}}
	    \put(39.375,167.5){\circle{4}}
	    
	    \put(78.75,10){\circle{4}}
	    \put(78.75,33){\circle{4}}
	    \put(78.75,88.75){\circle{4}}
	    \put(78.75,144.4){\circle{4}}
	    \put(78.75,167.5){\circle{4}}
	    
	    \put(118.125,10){\circle{4}}
	    \put(118.125,33){\circle{4}}
	    \put(118.125,88.75){\circle{4}}
	    \put(118.125,144.4){\circle{4}}
	    \put(118.125,167.5){\circle{4}}
	    
	    \put(157.5,10){\circle{4}}
	    \put(157.5,33){\circle{4}}
	    \put(157.5,88.75){\circle{4}}
	    \put(157.5,144.4){\circle{4}}
	    \put(157.5,167.5){\circle{4}}

	    % Labelling Drawing
	    \put(-3,10){\makebox(0,0)[r]{\small{$\phi_k^1$}}}
	    \put(-3,33){\makebox(0,0)[r]{\small{$\phi_k^2$}}}
	    \put(-3,167.5){\makebox(0,0)[r]{\small{$\phi_k^{N_{\eta}}$}}}
	    \put(160.5,167.5){\makebox(0,0)[l]{\small{$\phi_k^{N_{\eta}N_{\xi}}$}}}
	    \put(39.375,12){\makebox(0,0)[b]{\small{$\phi_k^{N_{\eta}+1}$}}}
	    \put(39.375,170.5){\makebox(0,0)[b]{\small{$\phi_k^{2N_{\eta}}$}}}

	  \end{picture}
	}	    	
    \end{minipage}
    \caption{$N_{\xi}$, $N_{\eta}$, are the number of points in the $\xi$ and $\eta$ direction respectively.}
    \label{fig:numindex}
 \end{figure}
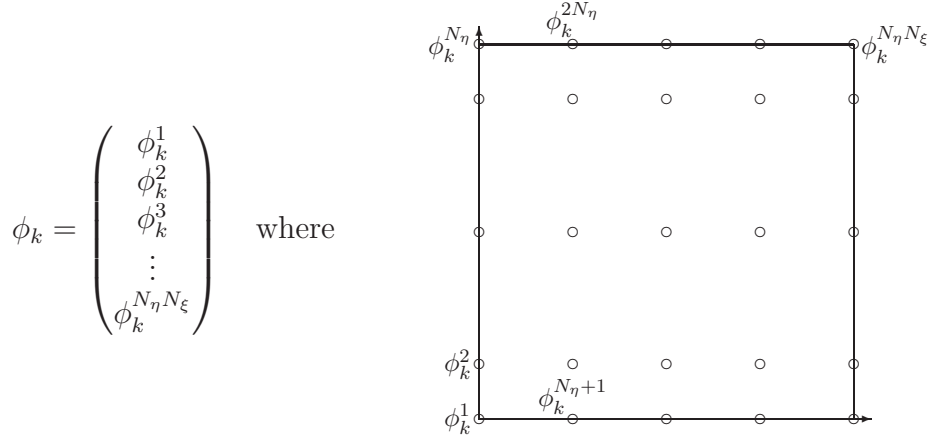

To solve (\ref{eq:symsystem}), (\ref{eq:numbound}) numerically we convert the variables $\phi_k$, $n_k$, $u_k^{\phi}$, $u_k^n$ and $G(\xi,\eta)$ from 2D continuous objects into vectors on the domain discretised with equispaced points in the periodic $\xi$ direction and Chebyshev spaced points in the non-periodic $\eta$ direction. To turn the 2D grid into a vector, we index from the bottom left hand corner as illustrated in figure \ref{fig:numindex}.

   To write the system (\ref{eq:symsystem}) as a matrix equation we introduce the notation:
\begin{align}
 D_{\xi}^n: & \quad \text{the $n$th order differentiation matrix in the $\xi$ direction}, \\
 \sd[n_k]:& \quad \text{a square matrix with the variable $n_k$ on its} \nonumber\\
 & \quad \text{diagonal.} 
\end{align}
We make use of the Schwarz-Christoffel toolbox \cite{driscoll2005algorithm} and take the differentiation matrices from the MATLAB dmsuite \cite{weideman2000matlab}, where $D_{\xi}^n$ is a Fourier differentiation matrix and $D_{\eta}^n$ is a Chebyshev differentiation matrix. Using this notation, (\ref{eq:symsystem}) becomes
{\footnotesize
\begin{align}
  \begin{pmatrix}
   \nu^2 \left( D_{\xi}^2 +  D_{\eta}^2\right) & \quad \sd[G(\xi,\eta)] \\
   \sd[D_{\xi}^1n_k]D_{\xi}^1 + \sd[D_{\eta}^1n_k]D_{\eta}^1 & D_{\xi}^2 + D_{\eta}^2 + \sd[D_{\xi}^1\phi_k]D_{\xi}^1 + \sd[D_{\eta}^1\phi_k]D_{\eta}^1 - \sd\left[\f{2n_k}{\nu^2}\cdot G(\xi,\eta)\right]
  \end{pmatrix}
  \begin{pmatrix}
   u_k^{\phi} \\
   u_k^n
  \end{pmatrix} \nonumber \\
  =-
  \begin{pmatrix}
   \nu^2 \left(D_{\xi}^2\phi_k + D_{\eta}^2\phi_k\right) + G(\xi,\eta) \cdot n_k, \\
   D_{\xi}^2n_k + D_{\eta}^2n_k + D_{\xi}^1n_k\cdot D_{\xi}^1\phi_k + D_{\eta}^1n_k\cdot D_{\eta}^1\phi_k - G(\xi,\eta)\cdot \f{n_k\cdot n_k}{\nu^2}
  \end{pmatrix}. \label{eq:numsystem}
\end{align}
}
where $\cdot$ denotes element-wise multiplication and the boundary conditions are
\begin{align}
  u^{\phi}(\xi,0) & = 0, & u^{\phi}(\xi,\eta^*) & = 0, & u^{n}(\xi,0) & = 0, & u^{n}(\xi,\eta^*) & = 0, \label{eq:boundnumf}
\end{align}
 Boundary conditions in the $\xi$ direction are not needed as we use Fourier differentiation matrices that assume periodicity.  To apply the boundary conditions (\ref{eq:boundnumf}) we replace rows in the matrices in (\ref{eq:numsystem}) so that the relevant elements of $u_k^{\phi}$ and $u_k^n$ are selected and set to zero.  We first solve the problem in $1D$ using \textsf{Chebfun} and then extend this into the second dimension to create the initial guess.  The system (\ref{eq:numsystem}) is then solved using damped Newton iteration where the damping parameter is calculated using the algorithm given in the appendix of \cite{Birkisson12}.

Finally, to obtain the variables $\phi$ and $n$ in the transformed frame we first interpolate the periodic points to Chebyshev points using fourint.m \cite{weideman2000matlab} and then use \textsf{Chebfun} to create \textsf{chebfun} objects of $\phi$ and $n$.  This enables us to access the wide range of functionality built into \textsf{Chebfun} for further analysis of the results. 

\subsection{Validation}

\begin{figure}[hbtp!]
\centering
  \includegraphics[width=0.9\textwidth]{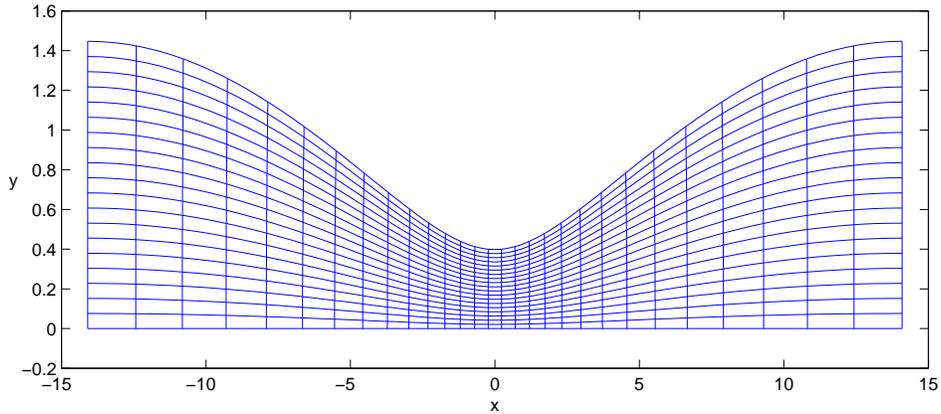}
  \caption{Domain used for validation.}
  \label{fig:validdomain}
\end{figure}

We now validate our numerical method by investigating its convergence for both the conformal map formulation and the hodograph plane formulation for the same test problem. The domain considered is produced using (\ref{eq:confmap}) with $L\approx 28.2$, $\epsilon=0.57$, and $\eta^*=0.84$ and is shown in figure \ref{fig:validdomain}. We also let $\Phi=1$ and $\nu=0.2$ in our numerical test. We are able to use the results from the conformal map formulation to determine the corresponding $QL$ and $y=F(\psi)$ in the hodograph plane formulation. 

To examine the convergence of the numerical method we fix the number of points in one direction and vary the number of points in the other and then observe how the integrated quantity $Q$ converges by considering the quantities  
\begin{align}
 \Delta Q^{N_{\xi}} & = |Q^{N_{\xi}+4} - Q^{N_{\xi}}|, & \Delta Q^{N_{\eta}} & = |Q^{N_{\eta}+4} - Q^{N_{\eta}}|,
 \end{align}
 for the conformal map approach, and
 \begin{align}
 \Delta Q^{N_v} & = |Q^{N_v+4} - Q^{N_v}|, & \Delta Q^{N_{\psi}} & = |Q^{N_{\psi}+4} - Q^{N_{\psi}}|.
\end{align}
for the hodograph approach.  

Tables \ref{tab:ConvergNxi}--\ref{tab:ConvergNpsi} show the results of our convergence tests and clearly demonstrate that both numerical methods rapidly converge to the same solution. We observe that the error decreases exponentially with the number of discretisation points, as expected for a spectral method.  It is evident that fewer points are required in the conformal map formulation than the hodograph plane formulation to achieve the same accuracy.  We believe that this occurs because the solutions for $\phi$ and $n$ are smoother in the conformal mapping plane than in the hodograph plane.  Additionally, we find in the hodograph plane formulation the numerical method is sensitive to the smoothness of the function $y(0,\psi)=F(\psi)$.

\begin{table}[hbtp!]
\footnotesize
\caption{Convergence Tests}
\centering
  \subtable[Convergence of conformal map formulation in $\xi$ direction, where $N_{\eta} = 31$.]{
 % \vspace{24pt}
  \centering
    \begin{tabular}{ | l | l | l |}
      \hline
      $N_{\xi}$ & $Q^{N_{\xi}}$ & $\Delta Q^{N_{\xi}}$ \\ \hline
      $4$& $0.63598$ & $8.5 \times 10^{-5}$ \\ \hline
      $8$ & $0.63607$ & $5.7 \times 10^{-8}$ \\ \hline
    \end{tabular}
   \label{tab:ConvergNxi}
   }
   \qquad
   \subtable[Convergence of conformal map formulation in $\eta$ direction, where $N_{\xi} = 12$.]{
    \centering
  \begin{tabular}{| l | l | l |}
    \hline
    $N_{\eta}$ & $Q^{N_{\eta}}$ & $\Delta Q^{N_{\eta}}$ \\ \hline
    $7$ & $0.46937$ & $1.6 \times 10^{-1}$ \\ \hline
    $11$ & $0.62815$ & $7.6 \times 10^{-3}$ \\ \hline
    $15$ & $0.63574$ & $3.2 \times 10^{-4}$ \\ \hline
    $19$ & $0.63605$ & $1.2 \times 10^{-5}$ \\ \hline
    $23$ & $0.63607$ & $4.4 \times 10^{-7}$ \\ \hline
    $27$ & $0.63607$ & $1.5 \times 10^{-8}$ \\ \hline
  \end{tabular}
  \label{tab:ConvergNeta}
  }
  \\% \vspace{12pt}
    \subtable[Convergence of hodograph plane formulation in $v$ direction, where $N_{\psi} = 32$.]{
  \centering
    \begin{tabular}{ | l | l | l |}
    \hline
    $N_v$ & $Q^{N_v}$ & $\Delta Q^{N_v}$ \\ \hline
    $17$& $0.63110$ & $4.5 \times 10^{-3}$ \\ \hline
    $21$& $0.63561$ & $4.0 \times 10^{-4}$ \\ \hline
    $25$ & $0.63600$ & $5.0 \times 10^{-5}$ \\ \hline
    $29$ & $0.63606$ & $7.0 \times 10^{-6}$ \\ \hline
    $33$ & $0.63606$ & $1.1 \times 10^{-6}$ \\ \hline
    $37$ & $0.63607$ & $2.1 \times 10^{-7}$ \\ \hline
    $41$ & $0.63607$ & $4.7 \times 10^{-8}$ \\ \hline
  \end{tabular}
   \label{tab:ConvergNv}
   }
   \qquad
   \subtable[Convergence of hodograph plane formulation in $\psi$ direction, where $N_v = 45$.]{
    \centering
  \begin{tabular}{ | l | l | l |}
    \hline
    $N_{\psi}$ & $Q^{N_{\psi}}$ & $\Delta Q^{N_{\psi}}$ \\ \hline
    $4$ & $0.57976$ & $5.1 \times 10^{-2}$ \\ \hline
    $8$ & $0.63071$ & $4.8 \times 10^{-3}$ \\ \hline
    $12$ & $0.63554$ & $4.7 \times 10^{-4}$ \\ \hline
    $16$ & $0.63601$ & $4.7 \times 10^{-5}$ \\ \hline
    $20$ & $0.63606$ & $4.8 \times 10^{-6}$ \\ \hline
    $24$ & $0.63607$ & $4.8 \times 10^{-7}$ \\ \hline
    $28$ & $0.63607$ & $4.9 \times 10^{-8}$ \\ \hline
  \end{tabular}
  \label{tab:ConvergNpsi}
  }
\end{table}

\subsection{Example numerical solutions}
\subsubsection{Smooth domain}
\label{sec:slownum}
\begin{figure}[hbtp!]
  \centering
    \includegraphics[width = 0.9\textwidth]{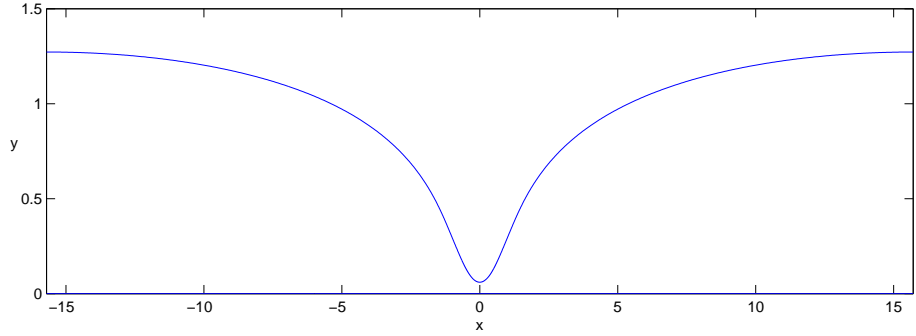}
    \caption{Example slowly varying domain used; $h_{\text{min}} = 0.06$ and the radius of curvature at the minimum is $a \approx 1.7$.}
    \label{fig:slowdomain}
\end{figure}
\begin{figure}
\centering
 \subfigure[]{
   \includegraphics[width=0.9\textwidth]{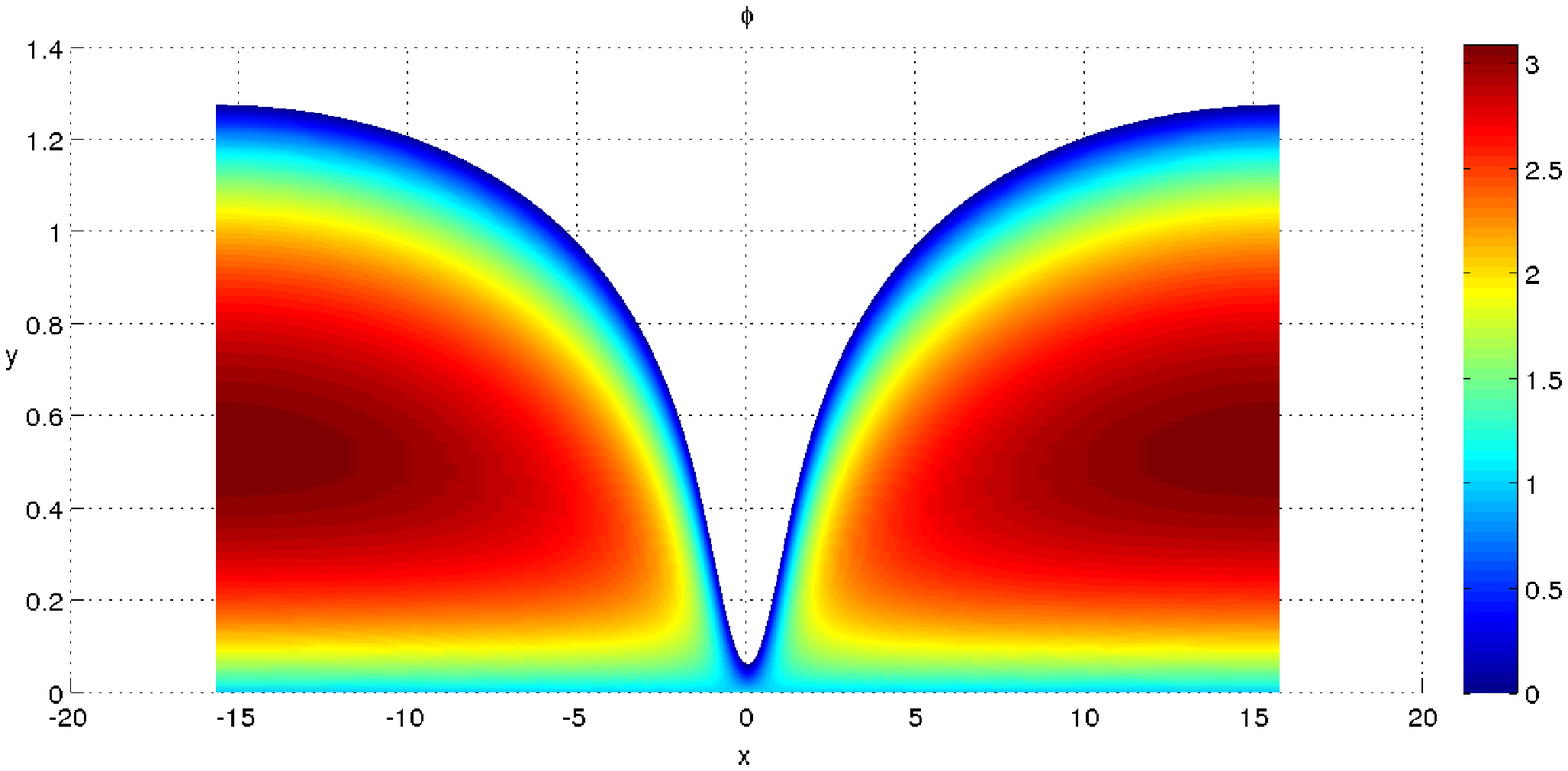}
   \label{fig:slowphi}
   }
  \\
  \subfigure[]{
   \includegraphics[width=0.9\textwidth]{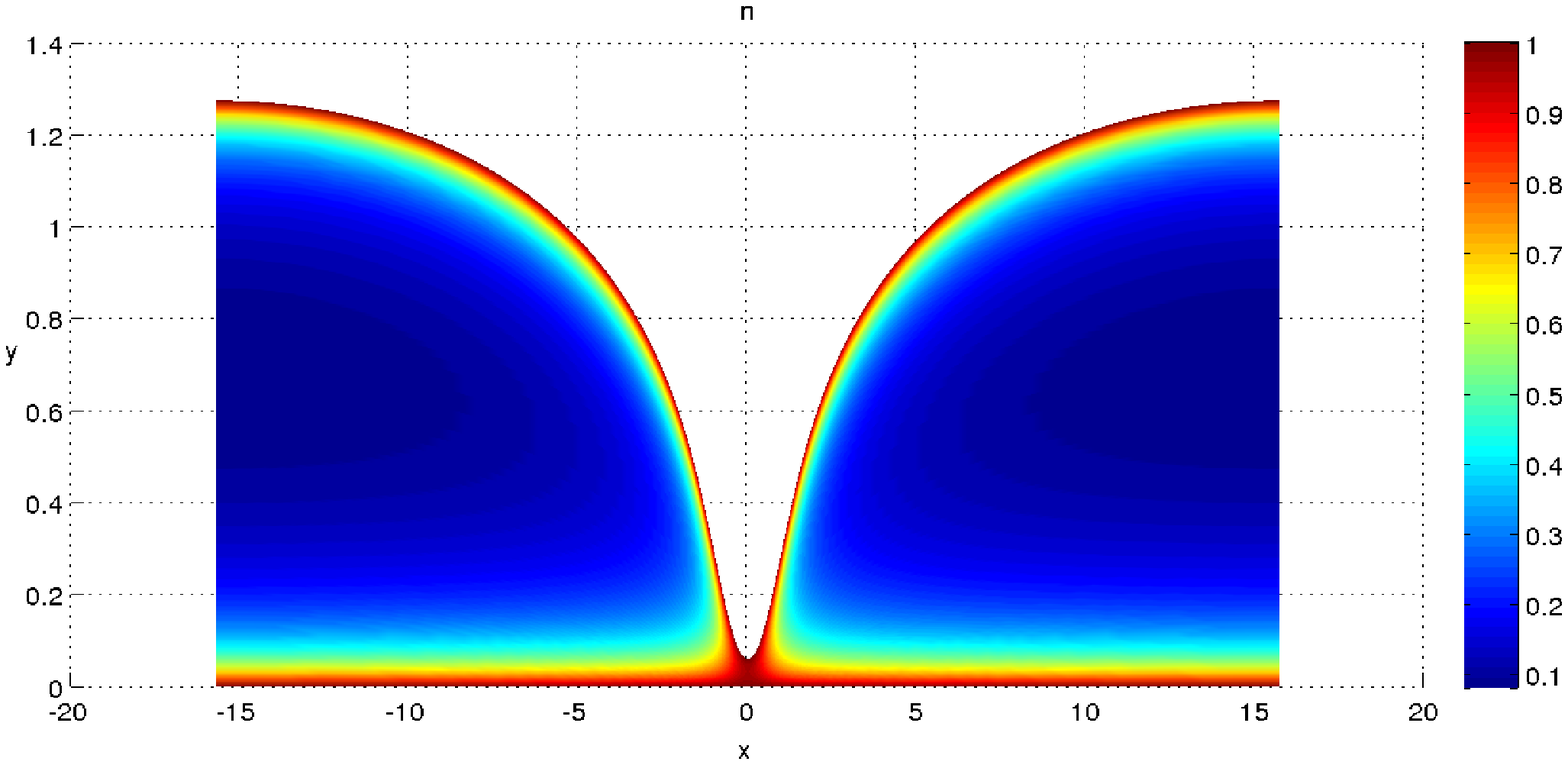}
   \label{fig:slown}
   }
  \caption{Electric potential $\phi$ and electron density $n$ plotted against position $x$ and $y$ with parameter values $\nu=0.1$, $\Phi=1$, $L\approx31.4$ and $Q=0.74$.}
  \label{fig:slowphiandn}
\end{figure}
\begin{figure}
\centering
  \subfigure[Full domain electron trajectories]{
    \includegraphics[width=0.9\textwidth]{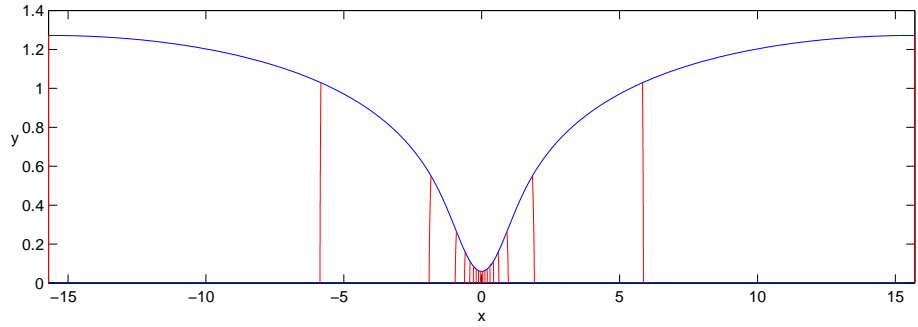}
    \label{fig:Trajectories1}
    }
  \\
  \subfigure[Central region electron trajectories]{
    \includegraphics[width=0.9\textwidth]{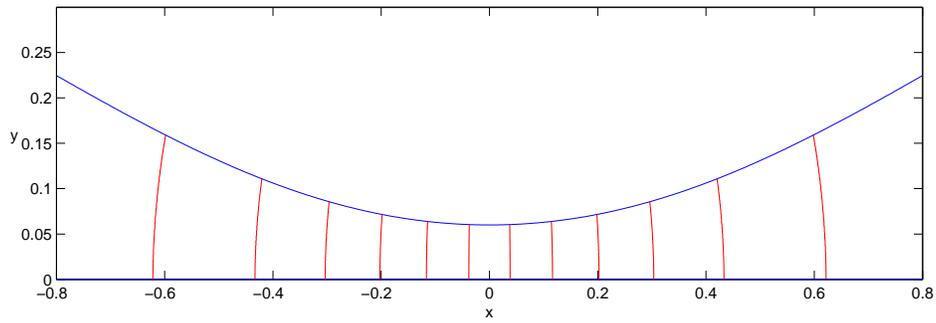}
    \label{fig:Trajectories2}
    }
  \caption{Electron trajectories through the glass layer.  An equal current is carried between each pair of adjacent trajectories.}
  \label{fig:Trajectories}
\end{figure}

We consider the domain shown in figure \ref{fig:slowdomain} where $h_{\text{min}}=0.06$, $L \approx 31.4$, the radius of curvature at $h_{\text{min}}$ is $a\approx 1.7$ and
\begin{equation}
F(\psi)=h_{\text{min}}\sec^2\left(\f{\psi}{\Phi}\sqrt{\f{h_{\text{min}}}{2a}}\right)
\end{equation}
in the central region, where the edges are patched to polynomials (this choice will be justified below in Section \ref{sec:smoothAsymp}).  For illustration, we set $\Phi=1$ and $\nu=0.1$, since we expect the normalised Debye length to be small in practice.  In figure \ref{fig:slowphiandn} we plot the numerical solutions for $\phi$ and $n$.  We observe boundary layers at the two surfaces of the glass layer away from the minimum $h=h_{\text{min}}$.  Near the minimum, however, we see that $n\approx 1$ while $\phi$ varies approximately linearly across the layer.  Both of these structures were found by Black et al. \cite{Black13} in limiting one-dimensional solutions.  The presence of the large electron density in the central region means we expect the majority of the current to be collected there.  This is visible in figure \ref{fig:Trajectories} where the electron trajectories are plotted with equal current between each trajectory and are largely concentrated around $h=h_{\text{min}}$.  
Also figure \ref{fig:Trajectories2} displays the quasi-1D nature of the solution near $h=h_{\text{min}}$, as the electron trajectories have little $x$-variation. This quasi-1D behaviour is due to the slowly varying top surface: note the different axes scalings in figure~\ref{fig:slowdomain}.  

The effective resistance of the domain in figure \ref{fig:slowdomain} with $\nu=0.1$ and $\Phi=1$ is $R = 1.35$.  On the otherhand, with $\nu=0.1$ and $\Phi=1$, a uniform  glass layer with the same average thickness has an effective resistance $R = 5.30$.  Therefore the short circuiting causes a reduction in the net resistance by a factor of approximately $4$.


\subsubsection{Wedge-like domain}

\begin{figure}[hbtp!]
\centering
  \includegraphics[width=0.9\textwidth]{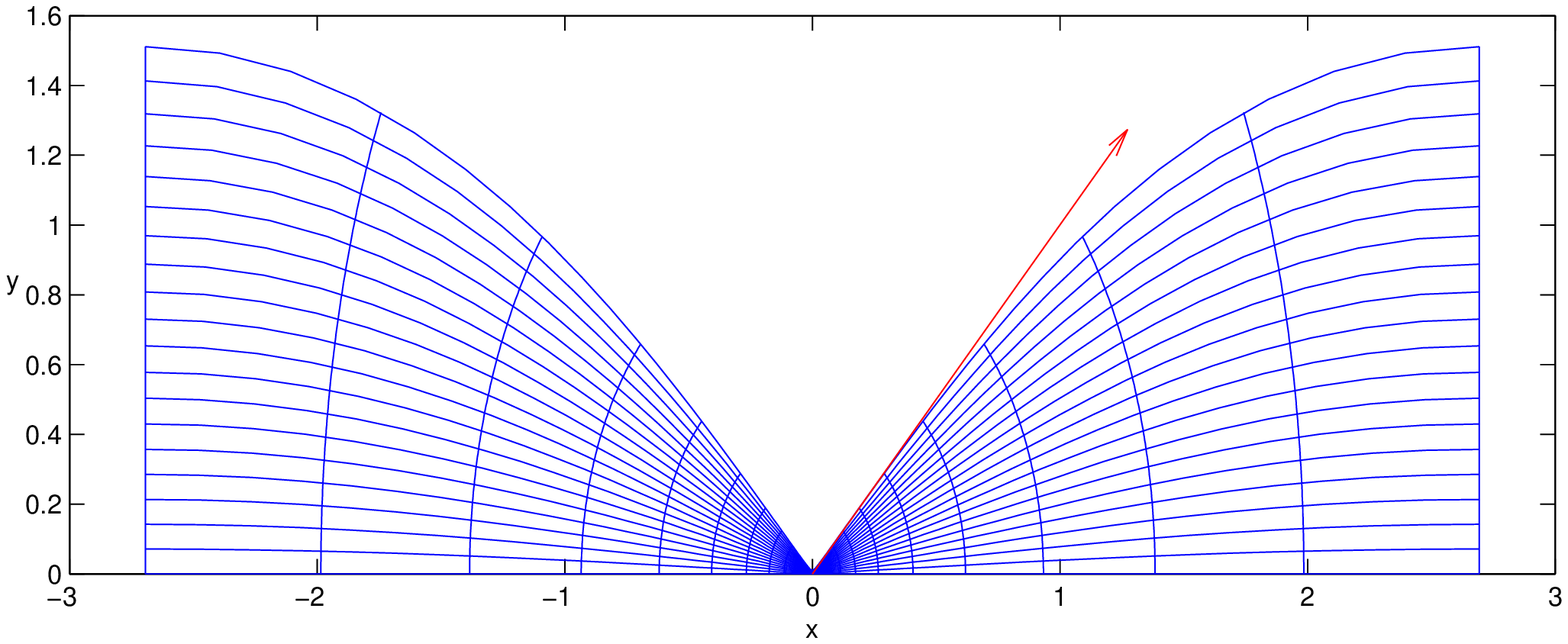}
  \caption{Example wedge domain used, $h_{\text{min}}=0.01$, $L\approx5.4$, $\epsilon \approx 0.016$, $\eta^*\approx 0.34$ and the arrow indicates the angle $\beta = \pi/4$.}
  \label{fig:wedgedomain}
\end{figure}
\begin{figure}[hbtp!]
\centering
  \includegraphics[width= 0.7\textwidth]{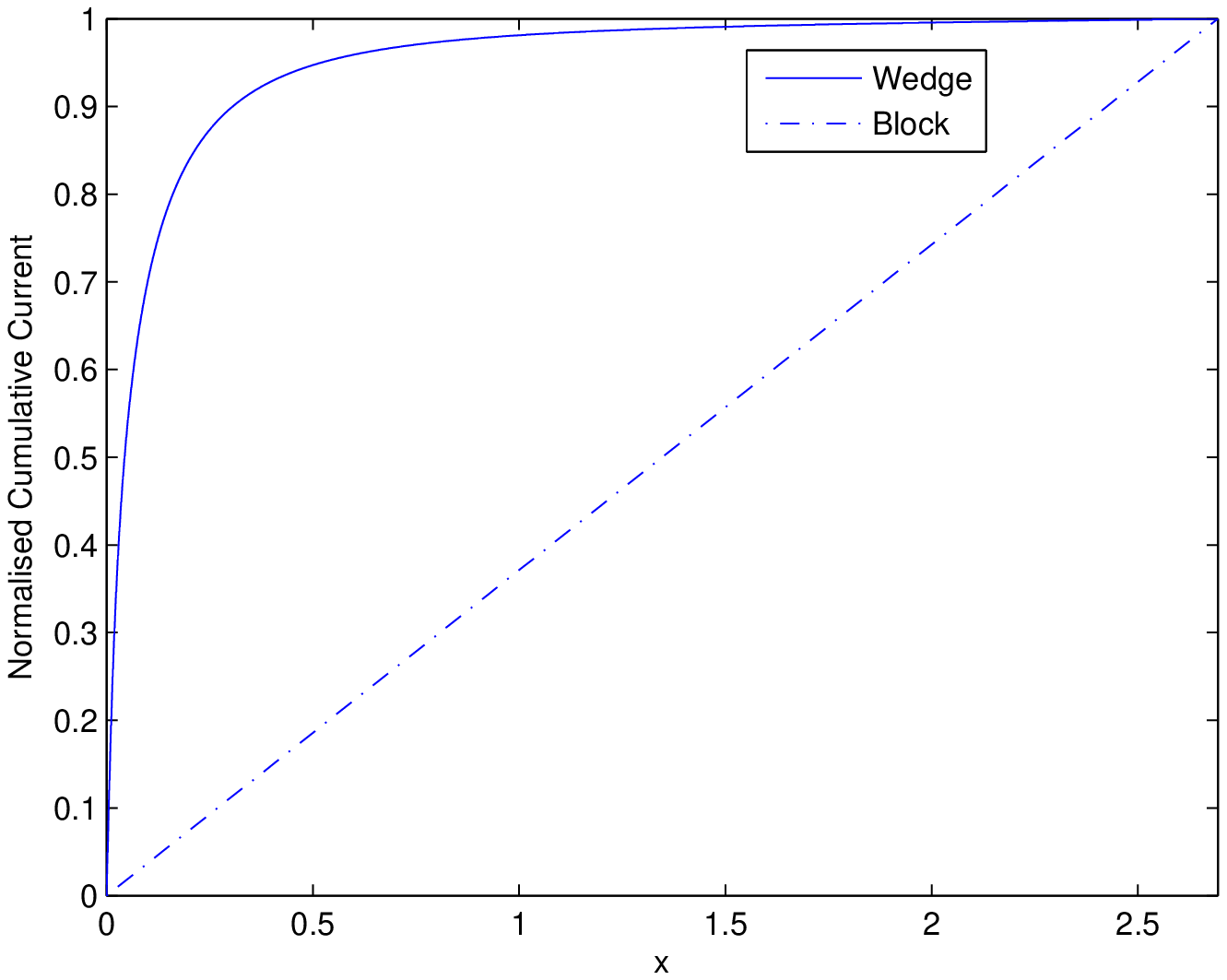}
  \caption{Normalised cumulative current, (\ref{eq:cumcurr}), plotted against horizontal position $x$. 'Wedge' is the normalised cumulative current for the domain shown in figure \ref{fig:wedgedomain}, and `Block' is the corresponding curve for a glass layer with constant thickness.}
  \label{fig:cumulative}
\end{figure}

We will now analyse a wedge domain where we let $\Phi=1$, $\nu =0.05$, $h_{\text{min}}=~0.01$, $L\approx5.4$ and $\beta = \pi/4$.  The corresponding parameter values in the mapping function (\ref{eq:confmap}) are $\epsilon \approx 0.016$ and $\eta^* \approx 0.34$. The domain is shown in figure \ref{fig:wedgedomain}.
\begin{figure}[hbtp!]
\centering
 \subfigure[]{
   \includegraphics[width=0.45\textwidth]{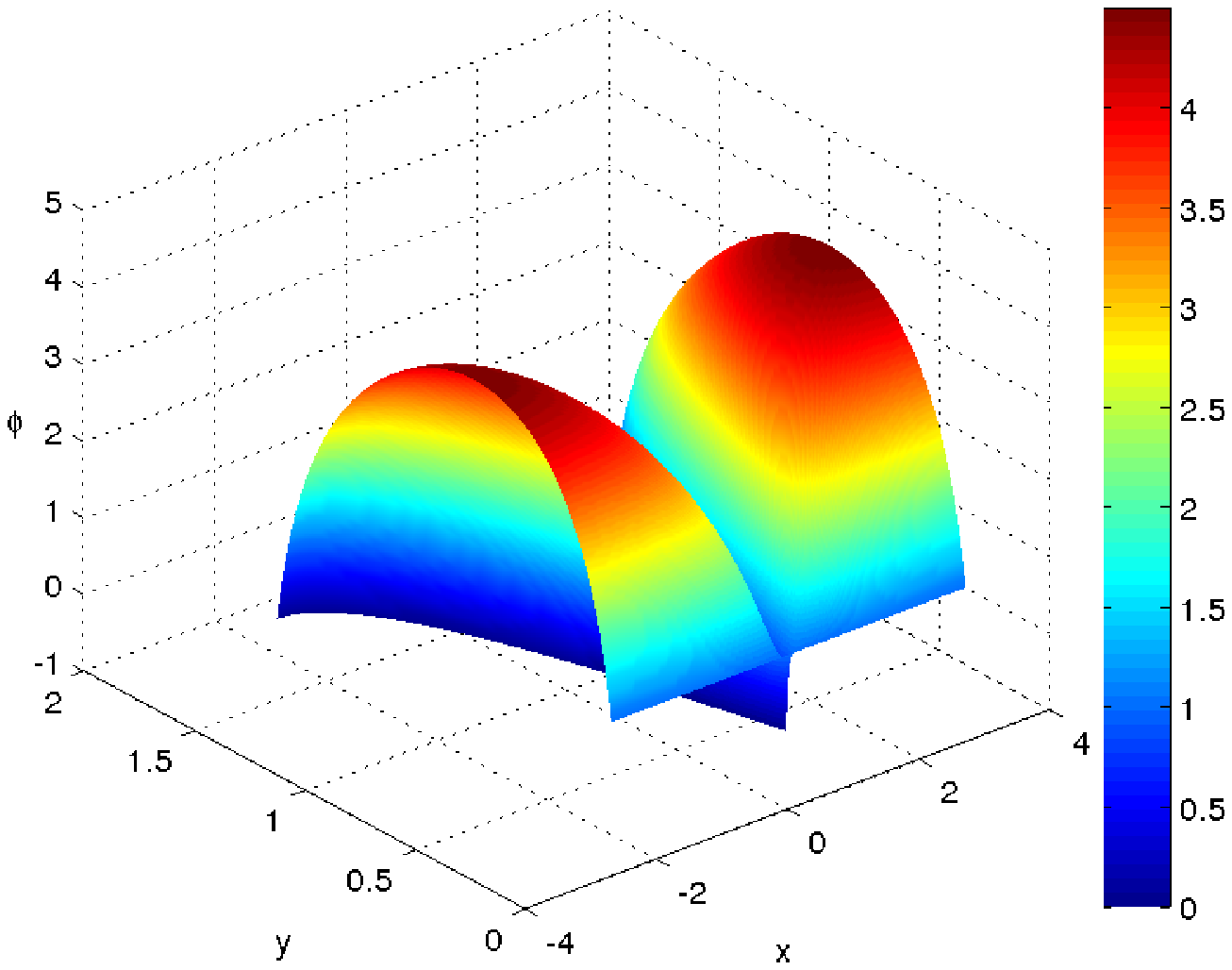}
   \label{fig:phiwedge}
   }
  \,
  \subfigure[]{
   \includegraphics[width=0.45\textwidth]{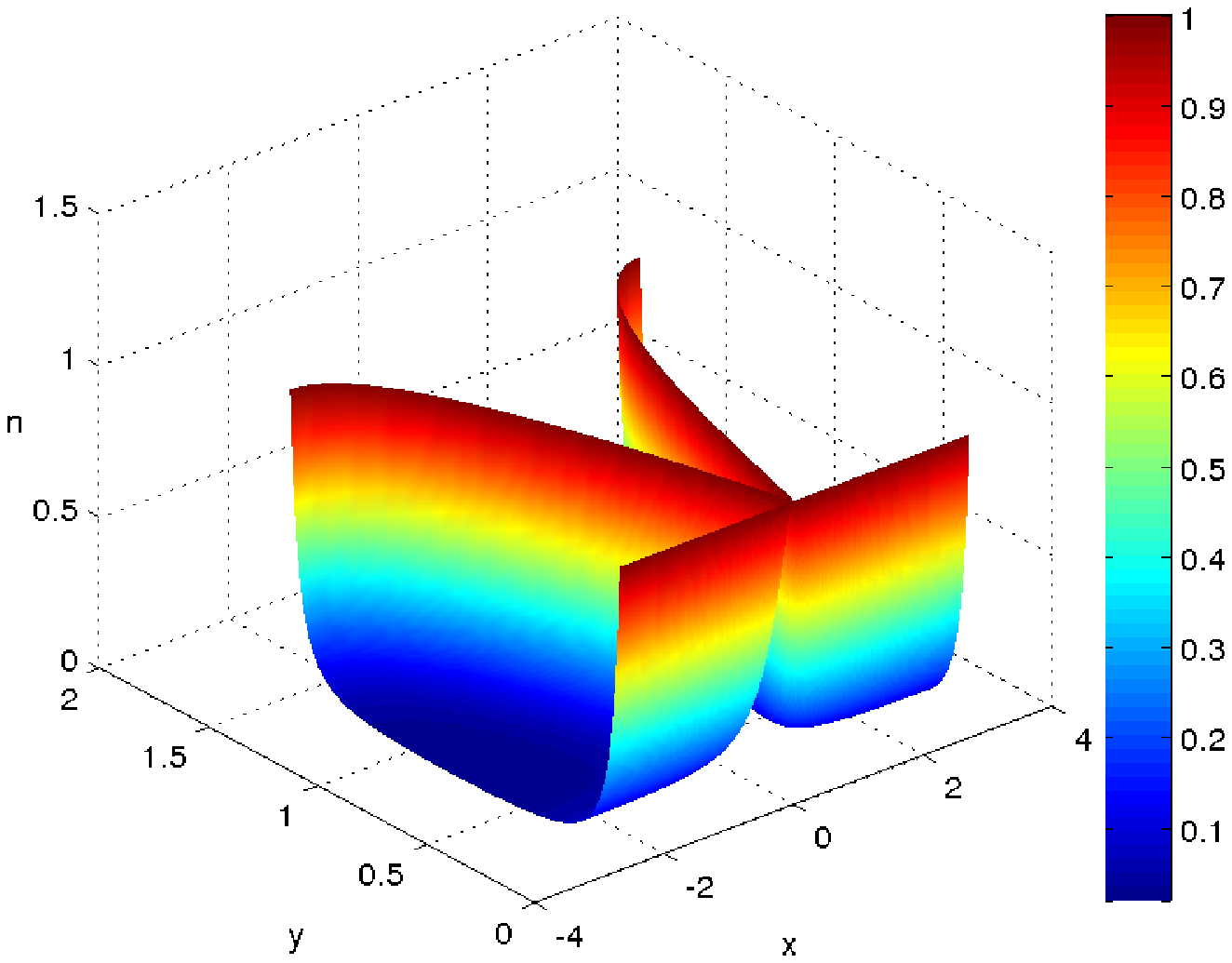}
   \label{fig:nwedge}
   }
  \caption{Electric potential $\phi$ and electron density $n$ plotted against position $x$ and $y$ with parameter values $\nu=0.05$, $\Phi=1$ and $Q=1.70$.}
  \label{fig:phiandnwedge}
\end{figure}
In figure \ref{fig:phiandnwedge} the numerical solutions for $\phi$ and $n$ are plotted.  These figures exhibit the same general features as figures \ref{fig:slowphi} and \ref{fig:slown}, with a boundary layer structure away from $h=h_{\text{min}}$ and a region of high current density near $h=h_{\text{min}}$. The short circuiting of the current through the thinner region is demonstrated by plotting the cumulative current distribution in figure \ref{fig:cumulative}. We observe that over $80\%$ of the current is collected within a neighbourhood $[-0.25,0.25]$ of $x=0$. In figure \ref{fig:cumulative} we include the corresponding curve for a glass layer with constant thickness for a reference.  In general, the more concave the normalised cumulative current is, the more the current is short-circuiting through the region around $h=h_{\text{min}}$.
Finally, the effective resistance for the domain shown in figure \ref{fig:wedgedomain}, with $\Phi=1$ and $\nu=0.05$, is $R = 0.6$.  For comparison, with the same values for $\Phi$ and $\nu$, the resistance of a uniform glass layer with the same average thickness is $R = 15.0$.  This example illustrates the dramatic effect that nonuniform thickness may have on the net resistance: here the short-circuiting of the current causes a reduction in the resistance by a factor of $25$.  The increased difference between these values when compared with those in Section \ref{sec:slownum} is mainly due to the smaller values of $h_{\text{min}}$ and $\nu$ used.

\subsection{Summary}

In this section we have demonstrated the implementation of a new numerical method that can be used to solve for the electron density and electric potential in both the conformal mapping and the hodograph plane formulations.  The method has been validated for both implementations, with exponential convergence evident.  The agreement that we find between the different implementations for the solution of the same problem gives high confidence in both methods.  We have considered some example cases and found that short-circuiting through the thinner regions of the glass layer causes a significant reduction in resistance when compared to a constant thickness domain.  This behaviour motivates us to seek asymptotic expressions for the average current density as the minimum thickness $h_{\text{min}}\rightarrow 0$.  We will see that the solution structure depends crucially on whether the local behaviour is ``slowly varying'' or ``wedge-like''.

\section{Asymptotic analysis}

\subsection{Slowly varying domain}
\label{sec:smoothAsymp}

We now calculate an approximation to the average current density in
the limit $h_{\text{min}} \rightarrow 0$ where the local radius of
curvature, $a\gg h_{\text{min}}$. It follows that the thickness profile near
the minimum is indeed slowly varying, with local aspect ratio
$\delta=(h_{\text{min}}/2a)^{1/2}\ll1$.
We focus on this region by performing the rescalings
\begin{align}
y&=h_{\text{min}}Y
&
x&=(2ah_{\text{min}})^{1/2}X.
\end{align}
Then the current density takes the form 
$\ve{j}=h_{\text{min}}^{-1}(\delta I,J)$
and, to lowest order in $\delta$, the problem
(\ref{eq:consvj})--(\ref{eq:bound}) is quasi-one-dimensional, i.e.
\begin{align}\label{quasi1Deq}
J_Y&=0,
&
J&=-n\phi_Y-n_Y,
&
(\nu/h_{\text{min}})^2\phi_{YY}&=-n,
\end{align}
subject to the boundary conditions 
\begin{align}\label{quasi1Dbc}
~&&
\phi&=\Phi,
&
n&=1 
& &
\text{at}\quad Y=0,
&&~
\nonumber\\
~&&
\phi&=0,
&
n&=1 
& &
\text{at}\quad Y=1+X^2. 
&&~
\end{align}

The problem now depends on $X$ only
parametrically, and the transformations
\begin{align}
\Yh&=\frac{Y}{1+X^2},
&
\Jh&=(1+X^2)J
&
\nuh&=\frac{\nu}{h_{\text{min}}(1+X^2)}
\end{align}
reduce (\ref{quasi1Deq})--(\ref{quasi1Dbc}) to the corresponding
purely one-dimensional problem
\begin{align}\label{1Deq}
\frac{\sd\Jh}{\sd\Yh}&=0,
&
\Jh&=-n\frac{\sd\phi}{\sd\Yh}-\frac{\sd n}{\sd\Yh},
&
\nuh^2\frac{\sd^2\phi}{\sd\Yh^2}&=-n,
\end{align}
\begin{align}\label{1Dbc}
\phi(0)&=\Phi,
&
n(0)&=1 
&
\phi(1)&=0,
&
n(1)&=1,
\end{align}
which was analysed in detail in \cite{Black13}. Let us denote the
solution for the current in this one-dimensional problem by
$\Jh=\jod\left(\nuh,\Phi\right)$.

By reversing the transformations carried out above, we deduce that the
net resistance of the layer in the limit $\delta\rightarrow0$ is
approximated by
\begin{equation}\label{eq:Qsmoothasymp2}
  Q \sim \f{1}{L} \sqrt{\f{2a}{h_{\text{min}}}}
\int_{-\infty}^{\infty}\jod\left(\f{\nu}{h_{\text{min}}(1 + X^2)},\Phi\right)
\,\frac{\sd X}{1 + X^2} .
\end{equation}
The function $\jod\left(\nuh,\Phi\right)$ is in general
determined numerically; however analytic 
approximations were found in various limiting cases in 
\cite{Black13}. In particular, in the limit $\nuh\rightarrow\infty$,
the layer acts as a resistor with $\jod\rightarrow\Phi$.
If $h_{\text{min}}\ll\nu$, i.e. the minimum thickness is smaller
than the Debye length, we may therefore use $\jod\approx\Phi$ to
approximate (\ref{eq:Qsmoothasymp2}) as
\begin{equation}\label{eq:slowavgcurrden}
 Q \sim \f{\pi\Phi}{L} \sqrt{\f{2a}{h_{\text{min}}}}.
\end{equation}

\begin{figure}[hbtp!]
\centering
 \includegraphics[width=0.9\textwidth]{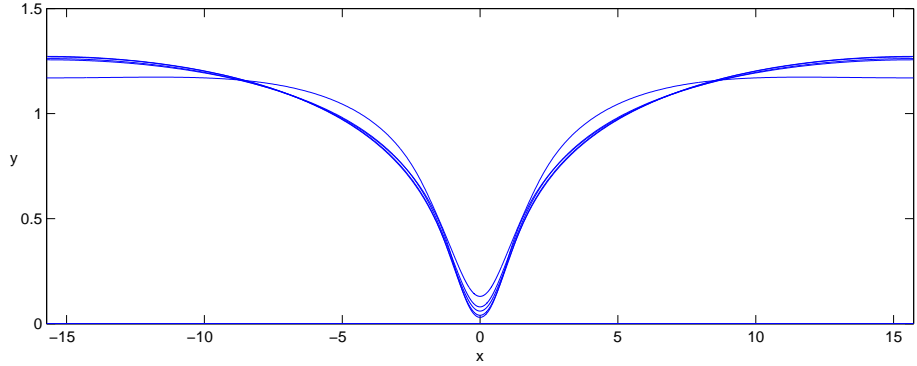}
 \caption{Slowly varying domains considered where $L\approx 31.4$, $h_{\text{min}} \in [0.03,0.13]$ and $a \in \approx~[1.6, 1.8]$.}
 \label{fig:slowdomains}
\end{figure}

Now we wish to test the approximation (\ref{eq:slowavgcurrden})
against our numerical results.
As noted in Section~\ref{sec:transsum}, the hodograph formulation is
most useful when \emph{a priori} knowledge of an appropriate $y=F(\psi)$ is known.
However, it requires us to pose the relation $y=F(\psi)$, where the
function $F(\psi)$ is indirectly related to the upper surface profile
$y=h(x)$. Here we can use the asymptotic behaviour of the solution
found above to determine the corresponding local behaviour of
$F(\psi)$ near the minimum thickness, namely by making use of
\begin{equation}
 \psi \sim \f{\Phi}{\delta} \int_0^X \f{\sd X}{1+ X^2},
\end{equation}
we determine
\begin{equation}
 F(\psi) \sim h_{\text{min}}
\sec^2\left(\f{\delta \psi}{\Phi}\right). \label{eq:yFpsi}
\end{equation}
Note that we expect $\psi=O(1/\delta)$ as $\delta\rightarrow0$.
Evidently (\ref{eq:yFpsi}) breaks down as
$\delta\psi\rightarrow\Phi\pi/2$ and $X\rightarrow\pm\infty$.
In any case, the numerical scheme requires $F(\psi)$ to be periodic,
which we achieve by patching the edges of (\ref{eq:yFpsi}) to
suitable polynomials.  A selection of domains used for comparison to our asymptotic expression is shown in figure \ref{fig:slowdomains}. As in the conformal map formulation the average thickness of the glass layer is always equal to $1$.

\begin{figure}[hbtp!]
\centering
  \includegraphics[width=0.7\textwidth]{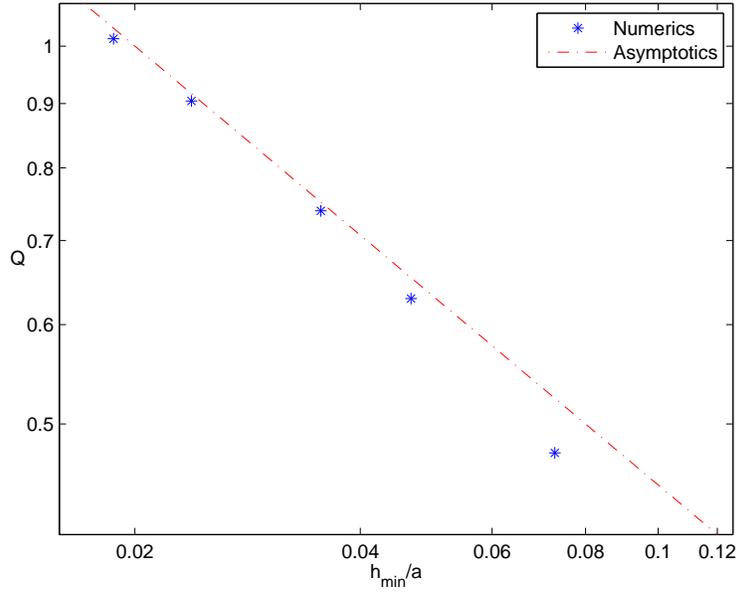}
  \caption{Average current density $Q$ plotted against glass layer
    minimum thickness $h_{\text{min}}/a$, where $\Phi=1$, $\nu=0.1$, $L \approx 31.4$ and $h_{\text{min}}/a \in \approx [0.019,0.073]$.
The asterisks show the results of numerical simulations; the
dot-dashed curve shows the asymptotic expression
(\ref{eq:slowavgcurrden}).}
  \label{fig:slowasymp}
\end{figure}

In figure~\ref{fig:slowasymp} we demonstrate that the asymptotic
expression (\ref{eq:slowavgcurrden}) gives an excellent approximation
to the results of numerical simulations for the average current
density $Q$ as $\delta\rightarrow0$.
This analysis accentuates how short-circuiting through any
thin spots in the glass layer may dominate the overall flux.
However, the approximation (\ref{eq:slowavgcurrden}) rests on the
assumption that $h_{\text{min}}\ll a$, and will therefore fail if the
curvature of the glass surface is too large in the neighbourhood of
the minimum. To describe the local behaviour in such sharp
``wedge-like'' domains, a different asymptotic approach is required,
as will be demonstrated below.

\subsection{Wedge-like domain}

\subsubsection{Introduction}

\begin{figure}[hbtp]
\centering
  \subfigure[Outer problem]{
   \centering
   \scalebox{0.8}{
      \begin{picture}(150,100)
%       % Drawing outline
	\put(0,15){\vector(1,0){140}}
	\put(70,15){\line(1,1){55}}
	\put(70,15){\line(-1,1){55}}
	\put(70,5){\vector(0,1){90}}
	\qbezier(100,15)(97.5,27.5)(91,36)
	\qbezier(125,70)(130,75)(135,74)
	\qbezier(15,70)(10,75)(5,74)
	
	% Labelling Drawing
	\put(82.5,22.5){\makebox(0,0)[l]{$\beta$}}
	\put(50,65){\makebox(0,0)[c]{$\phi=0$}}
	\put(90,65){\makebox(0,0)[c]{$n=1$}}
	\put(90,5){\makebox(0,0)[c]{$n=1$}}
	\put(50,5){\makebox(0,0)[c]{$\phi=\Phi$}}
	\put(70,97){\makebox(0,0)[b]{$y$}}
	\put(142,15){\makebox(0,0)[l]{$x$}}
	\put(0,5){\makebox(0,0)[c]{$-L/2$}}
	\put(135,5){\makebox(0,0)[c]{$L/2$}}
      \end{picture}
      }
     }
      \qquad 
  \subfigure[Inner problem]{
  \centering
  \scalebox{0.8}{
    \begin{picture}(150,100)
%       % Drawing outline
	\put(0,15){\vector(1,0){140}}	
	\qbezier(70,34)(90,34)(125.5,80)
	\qbezier(70,34)(50,34)(14.5,80)
	%\curve(12.5,85,70,34,127.5,85)
	\put(70,5){\vector(0,1){90}}
	\multiput(70,15)(5,5.4){13}{\line(0,1){1}}
	\multiput(70,13)(-5,5.4){13}{\line(0,1){1}}
	\arc(85,15){8}

	% Labelling Drawing
	\put(70,97){\makebox(0,0)[b]{$\tilde{y}$}}
	\put(142,15){\makebox(0,0)[l]{$\tilde{x}$}}
	\put(90,5){\makebox(0,0)[c]{$n=1$}}
	\put(50,5){\makebox(0,0)[c]{$\phi=\Phi$}}
	\put(50,65){\makebox(0,0)[c]{$\phi=0$}}
	\put(90,65){\makebox(0,0)[c]{$n=1$}}
	\put(93,23){\makebox(0,0)[c]{$\beta$}}

      \end{picture}
      }
   }
      \caption{Schematics of outer and inner problems.}
      \label{fig:wedgeasymps}
   \end{figure}
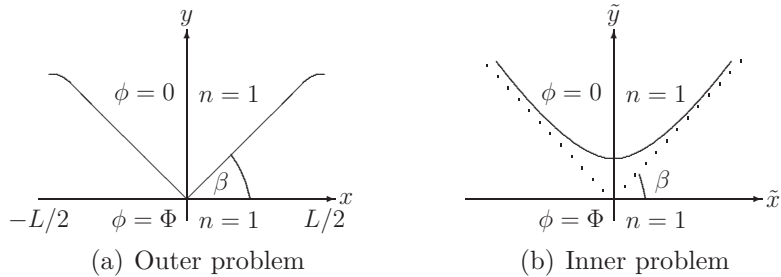

In this section we investigate the behaviour of the average current
density in a wedge-like domain, with $a = O(h_{\text{min}})$
as $h_{\text{min}} \rightarrow 0$.
If we naively set $h_{\text{min}}=0$, we are left with the ``outer''
problem illustrated in figure~\ref{fig:wedgeasymps}(a): the glass
layer thickness reaches zero with a corner singularity at (without
loss of generality) $x=0$.
The problem is regularised over a small region in which
$(x,y)=h_{\text{min}}\left(\tilde{x},\tilde{y}\right)$, as illustrated
in figure~\ref{fig:wedgeasymps}(b). In this ``inner'' problem, the
glass layer has unit minimum thickness and approaches a wedge shape as
$\left(\tilde{x},\tilde{y}\right)\rightarrow\infty$.
Below we will analyse the inner and outer problems and match them
asymptotically to obtain an approximation for the net flux through the
layer.

\subsubsection{Outer problem}
\label{sec:wedgeasympout}

The outer problem, namely the full governing equations
(\ref{eq:consvj})--(\ref{eq:jpois})
with the boundary conditions shown in figure~\ref{fig:wedgeasymps}(a),
requires numerical solution in general. For matching purposes, we only
require the asymptotic behaviour of the solution as $r\rightarrow0$,
which is easily found to be
\begin{align}\label{eq:wedgeoutphi}
 \phi & \sim
\Phi\left(1 - \f{\theta}{\beta}\right)+O\left(r^2\right), 
&
n & \sim 1+O\left(r^2\right),
\end{align}
where $(r,\theta)$ are plane polar coordinates.

\subsubsection{Inner problem}
\label{sec:wedgasympinn}

In the inner problem we scale $(x,y)=h_{\text{min}}\left(\xt,\yt\right)$ 
and let $h_{\text{min}} \rightarrow 0$ while expanding the solution as
\begin{align}
 n & \sim \nt_0 + h_{\text{min}}^2\nt_1 + \cdots,
&
 \phi & \sim \phit_0 + h_{\text{min}}^2\phit_1 + \cdots.
\end{align}
We find that $\nt_0=1$ and $\phit_0$ satisfies Laplace's equation
\begin{equation}
\frac{\partial^2\phit_0}{\partial\xt^2}
+\frac{\partial^2\phit_0}{\partial\yt^2}=0,
\end{equation}
subject to the boundary conditions
shown in figure~\ref{fig:wedgeasymps}(b).
In principle, this problem is easily solved by conformal mapping.
Assuming that the inner geometry shown in
figure~\ref{fig:wedgeasymps}(b) is the image of a strip
$0\leq\im\bigl(\zetat\bigr)\leq1$ under the conformal map
\begin{equation}
\zt=\xt+\ii\yt=\ft\bigl(\zetat\bigr)=\ft\bigl(\xit+\ii\etat\bigr),
\end{equation}
then the solution in the strip is given by
\begin{equation}\label{eq:wedgeinnasympphi}
\phit_0=\Phi\bigl(1-\etat\bigr).
\end{equation}
Here we have assumed that $\ft$ maps the real line to itself and also
fixes the points $0$, $\ii$ and $\infty$.
For example, the locally hyperbolic upper surface described by 
(\ref{eq:confmapinn})--(\ref{hyperbola})
corresponds to the mapping function
\begin{equation}\label{hyperft}
\ft\bigl(\zetat\bigr)=
\frac{\sinh\bigl(\beta\zetat\bigr)}{\sin\beta}.
\end{equation}

\subsubsection{Average current density approximation}

Now we find the leading-order average current density through the layer
by integrating over the inner and outer regions. Assuming symmetry
about the $y$-axis, we have
\begin{align}\label{eq:I1andI2}
 Q & = \f{2}{L} \left(\underbrace{
\int_{\lambda}^{L/2} \left[-n\frac{\partial\phi}{\partial y}
-\frac{\partial n}{\partial y}\right]_{y=0} \!\!\!\! \sd x}_{I_{out}}
+ \underbrace{\int_{0}^{\lambda/h_{\text{min}}}
\left[-\nt\frac{\partial\phit}{\partial\yt}
-\frac{\partial\nt}{\partial\yt}\right]_{\yt=0} \!\!\!\! \sd\xt}_{I_{in}} \right),
\end{align}
where $\lambda$ is assumed to satisfy $h_{\text{min}}\ll\lambda\ll 1$, and $I_{out}$, $I_{in}$ are the contributions from the outer and inner regions respectively.  We know from (\ref{eq:wedgeoutphi}) that the outer integral
$I_{\text{out}}$ diverges logarithmically as $\lambda\rightarrow0$
and, by subtracting the singular part, we obtain
\begin{equation}\label{Iout}
I_{\text{out}}\sim\int_0^{L/2} \left[-n\frac{\partial\phi}{\partial y}
-\frac{\partial n}{\partial y}-\frac{\Phi}{\beta x}\right]_{y=0}
\, \sd x+\frac{\Phi}{\beta}\log\left(\frac{L}{2\lambda}\right).
\end{equation}

For the inner contribution, we may perform the integral in the
$\zetat$-plane to obtain
\begin{equation}
I_{\text{in}}=\int_0^{{\xit^*}}
\left[-\nt\frac{\partial\phit}{\partial\etat}
-\frac{\partial\nt}{\partial\etat}\right]_{\etat=0} \!\!\!\! \sd\xit
\sim\Phi\xit^*,
\end{equation}
where $\xit^*$ is the point on the $\xit$-axis corresponding to
$\xt=\lambda/h_{\text{min}}$, i.e.
\begin{equation}
\ft\bigl(\xit^*\bigr)=\frac{\lambda}{h_{\text{min}}}.
\end{equation}
Now, from the wedge geometry, we must have
$\ft\bigl(\zetat\bigr)\sim b\mathrm{e}^{\beta\zetat}$ as
$\zetat\rightarrow\infty$, where $b$ is a real constant that depends
on the detailed inner geometry;
for example, $b=1/(2\sin\beta)$ for the mapping function  
(\ref{hyperft}).
Hence the leading-order inner integral is given by
\begin{equation}\label{Iin}
I_{\text{in}}\sim \frac{\Phi}{\beta}
\log\left(\frac{\lambda}{bh_{\text{min}}}\right).
\end{equation}

Finally, by substituting (\ref{Iout}) and (\ref{Iin}) into
(\ref{eq:I1andI2}), we obtain the following approximation for the
average flux:
\begin{equation}
Q\sim\frac{2\Phi}{\beta L}\log\left(\frac{L}{h_{\text{min}}}\right)
+C, \label{eq:wedgeavgQ}
\end{equation}
where the order-one constant $C$ is given by
\begin{equation}\label{eq:wedgeasympcurrden}
C=\frac{2}{L}\int_0^L \left[-n\frac{\partial\phi}{\partial y}
-\frac{\partial n}{\partial y}-\frac{\Phi}{\beta x}\right]_{y=0}
\, \sd x-\frac{2\Phi}{\beta L}\log\left(2b\right).
\end{equation}

\subsubsection{Comparison between asymptotics and numerics}
\begin{figure}[hbtp!]
 \centering
 \includegraphics[width=0.7\textwidth]{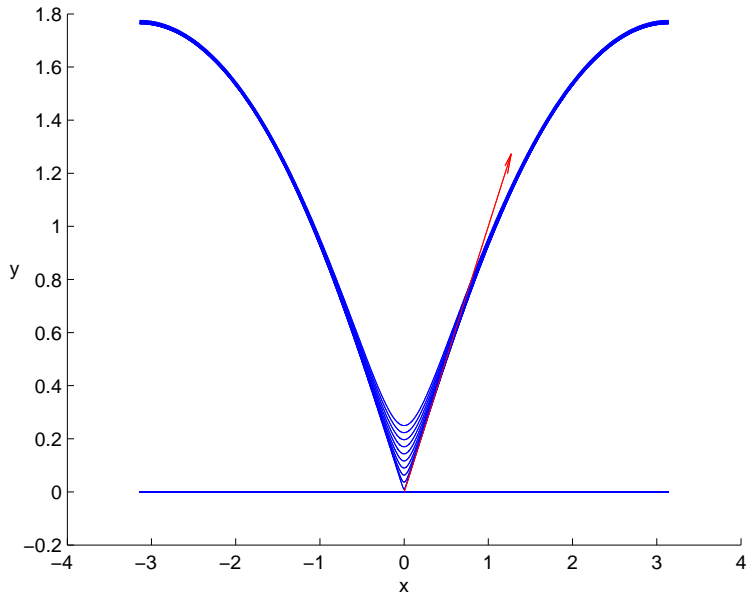}
 \caption{Wedge domains considered where $h_{\text{min}} \in [0.01, 0.25]$, $L = 2\pi$ and the arrow indicates the angle $\beta = \pi/4$.}
 \label{fig:wedgedomains}
\end{figure}

\begin{figure}[hbtp!]
 \centering
  \includegraphics[width=0.7\textwidth]{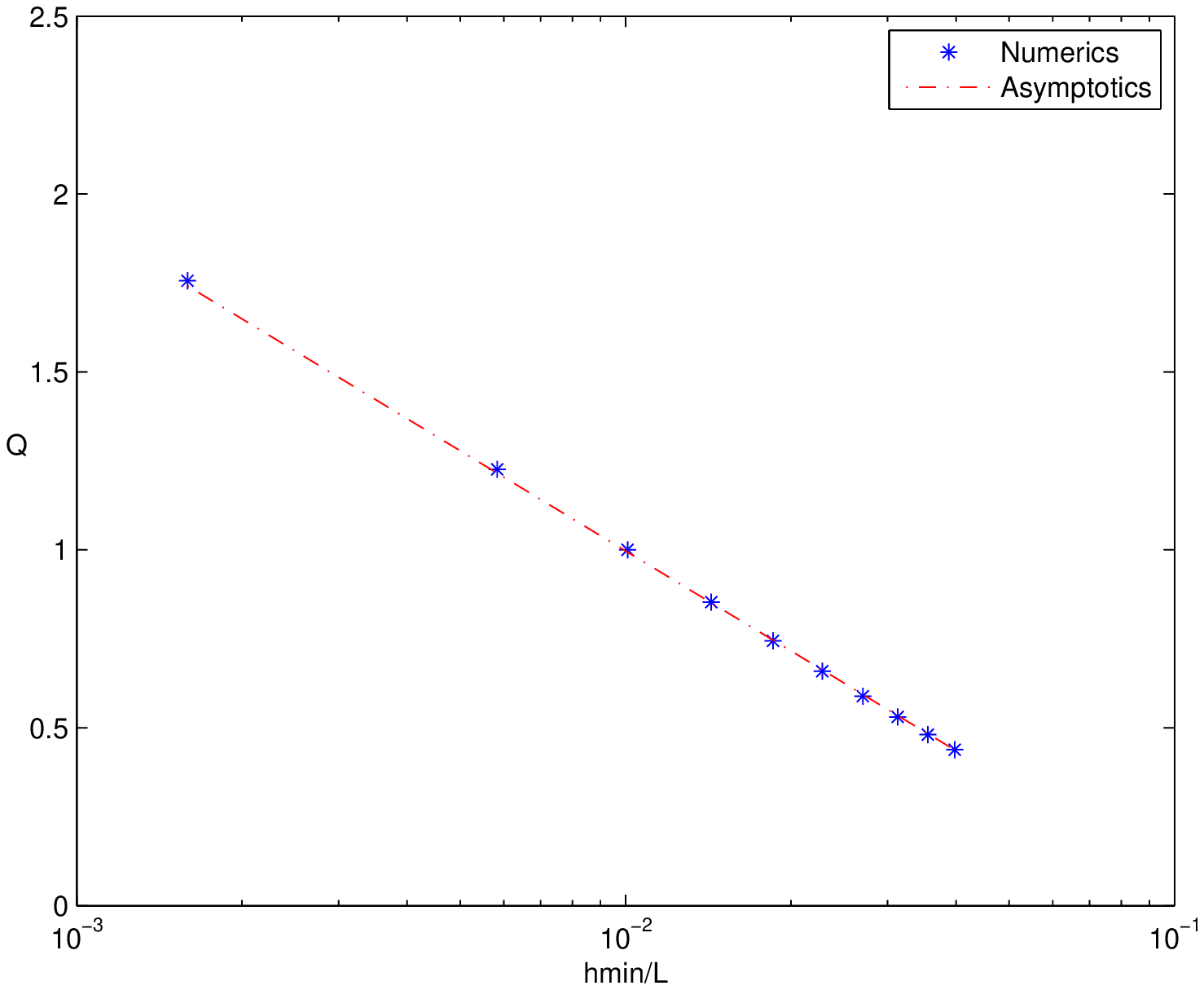}
  \caption{Average current density $Q$ plotted against glass layer minimum thickness $h_{\text{min}}$, where $\beta=\pi/4$, $\Phi=1$, $\nu = 0.1$, $L  = 2\pi$ and $h_{\text{min}} \in [0.01,0.25]$. The asterisks show the results of numerical simulations; the dot-dashed curve shows the asymptotic expression (\ref{eq:wedgeavgQ}) where we fit to find for this geometry $C \approx -0.87$ .}
  \label{fig:wedgeavgcurr}
\end{figure}

To validate our asymptotics we perform a sequence of simulations that fix the constant $C$.  We therefore consider a range of domains with the same $\beta$ and $L$, as shown in figure \ref{fig:wedgedomains}, and take $\nu=0.1$ and $\Phi=1$ in each simulation.  We must then relax the condition for the average thickness to be $1$.  The results are plotted in the linear-log plot in figure \ref{fig:wedgeavgcurr} where it is evident  the asymptotic expression (\ref{eq:wedgeavgQ}) gives a very good approximation for the average current density $Q$ as $h_{min}/L \rightarrow 0$.  In principle, one could calculate $C$ by solving the outer problem but we curve-fit to find $C \approx -0.87$ for this particular geometry.  Similarly to Section (\ref{sec:smoothAsymp}) the analysis demonstrates how the overall flux is dominated by short circuiting through thinner regions.

\section{Conclusion and discussion}

In this paper we extend to two dimensions the mathematical model from
\cite{Black13} for steady electron transport through a glass layer
between two electrodes.
We maintain many of the same modelling
assumptions: that the charge is predominantly carried by electrons,
that the electron flow is governed by drift and diffusion,
and that the electron
densities are known either side of the glass layer.  The analysis
focuses on the short-circuiting of current through thinner regions of
the glass layer, and we present a spectral numerical method to solve
the model once it has been mapped onto a rectangular domain.

We make use of two different mapping techniques;
which approach is preferable depends on the geometry under
consideration and the questions one wishes to answer.
In both cases, we validate the numerical method through its
exponentially rapid
convergence and its ability to produce insightful
results such as the electron trajectories and the normalised
cumulative current.
In considering nonuniform domains, the numerical results demonstrate
how the current short-circuits through thinner regions of the
glass. Consequently, the effective resistance is significantly lower
than comparable domains with constant glass layer thickness.

We use the numerical simulations to inform asymptotic calculations
conducted on domains with small minimum thickness $h_{\text{min}}$.
We consider two canonical local geometries.
First we suppose that the thickness profile is smooth, and the
local layer thickness is therefore ``slowly varying'' as
$h_{\text{min}}\rightarrow0$. In this limit, we find that the current
through the layer diverges like
$Q\propto\left(a/h_{\text{min}}\right)^{1/2}$,
where $a\gg h_{\text{min}}$ is the local radius of curvature.
The second regime considered is where $a=O\left(h_{\text{min}}\right)$,
so that the local behaviour of the
surface is ``wedge-like''. In this latter case, we find a weaker
logarithmic divergence of the current as $h_{\text{min}}\rightarrow0$.
Both predictions are found to agree well with full numerical
calculations.

The generalisation of our model into two dimensions is
well warranted. In experimental images of the front contact (see for
example \cite{Li09, Li11}), the glass layer is observed to vary
greatly in thickness, from $\sim 10$nm--$1\mu$m.
However, some of our modelling assumptions should be readdressed, such
as the neglect of holes and the simplified boundary conditions
(discussed in more detail in \cite{Black13}).  The presence of a positive charged species would require another transport law and the addition of a thermodynamic equilibrium law; an example is given in \cite{Please82}.  A physical process we have not considered is that silver ions could intercalate and transport across the glass layer leading to a dendritic short circuit.
Additionally, we note that, in the limit as the glass layer becomes very thin,
the validity of the continuum model must eventually be questioned.
Preliminary work suggests that
below $10$nm quantum tunnelling effects start to play a role in
determining the resistance of the glass layer.  These quantum effects
can be modelled through a modified drift-diffusion model, with a
higher order derivative term added to (\ref{eq:jdrifdif}); see, for
example \cite{Ancona87, cumberbatch2006nano, Pinnau02}.
Finally, it is well known that at
very high electric fields the drift-diffusion equations lose their
validity \cite{sze2006physics}.
As the glass layer thickness decreases below
$10$nm, sufficiently high fields could be present at relatively
moderate potential differences.
In these high field limits it is possible instead to use Monte
Carlo methods to simulate charge transport; see
\cite{jacoboni1989monte}.

Our model could be the building block for studying other, related, physical problems with some modification.  For example, in electrochemical thin films where there are multiple charged species moving under drift and diffusion, more complicated boundary conditions such as Butler-Volmer would be needed to describe reaction kinetics for electron transfer \cite{Bazant05current, Bonnefort01, Chu2005electrochemical}.

The dominant conduction mechanism for electron transport across the glass layer is under fierce debate and mathematical models such as the one explored in this paper can help with hypothesis testing.  In particular, our study has shown that the geometry of the glass layer makes significant difference to electron transport and hence performance of the photovoltaic cell.
\section*{Acknowledgements}

The authors are grateful to \'{A}. Birkisson, N. Hale, C. P. Please and the employees of DuPont (UK) Ltd for many useful discussions.  This work is supported by EPSRC and DuPont (UK) Ltd.  

\bibliographystyle{plain}	
\bibliography{ref}
\end{document}